\documentclass[10pt]{article} 

\usepackage{amsmath}
\usepackage{amssymb} 
\usepackage{latexsym} 
\usepackage{epic,eepic,epsfig,psfrag} 
\usepackage{amscd}  
\usepackage{index}

\newcommand{\QED}{\hspace*{\fill}$\Box$\medskip} 
\def\one{\hbox{1\hskip-2.7pt l}}
\def\smone{{\scriptstyle\rm 1\hskip-2.05pt l}}

\def\rd{{\rm d}}
\def\rT{{\rm T}}

\def\rG{{\rm G}}
\def\p{\phi}

\def\a{\alpha}
\def\b{\beta} 
\def\d{\delta} 

\def\ep{\varepsilon} 
\def\e{\eta} 

\def\th{\theta} 
\def\r{\rho}
\def\s{\sigma} 
\def\t{\tau} 
 
\def\x{\xi} 
 
\def\n{\nu} 
 
\def\z{\zeta} 
\def\o{\omega} 
 
\def\D{\Delta} 
 
\def\G{\Gamma} 
\def\pg{\mathhexbox278}
\def\S{\Sigma}
\def\Si{\Sigma} 
\def\Th{\Theta} 
 
\def\Om{\Omega} 
\def\P{\Phi} 
\def\cg{\mathfrak{g}}
\def\cA{{\mathcal A}}
 
\def\cC{{\mathcal C}}

\def\cG{{\mathcal G}}

\def\cL{{\mathcal L}}

\def\cU{{\mathcal U}} 
\def\cV{{\mathcal V}}
\def\cW{{\mathcal W}}

\def\R{{\mathbb R}}

\def\N{{\mathbb N}} 
 
\def\H{{\mathbb H}} 
\def\Z{{\mathbb Z}}
 
\def\half{{\textstyle{\frac 12}}} 
\def\im{{\rm im}\,}
\def\laplace{\Delta} 
\def\Pr{{\bf Proof:}\;} 
\def\st{\: \big| \:} 
\DeclareMathOperator{\supp}{supp}

\def\dt{{\rm d}t}
\def\ds{{\rm d}s}
\def\pd{\partial}

\def\comp{\circ}

\def\tp{{\tilde p}}

\def\tA{{\tilde A}}

\def\tU{{\tilde U}}
\def\la{\langle\,}
\def\ra{\,\rangle}

\def\HF{{\rm HF}}

\def\SU{{\rm SU}}

\def\tsum{\textstyle\sum}
\def\tint{\textstyle\int}

\newtheorem{dfn}{Definition}[section] 
\newtheorem{lem}[dfn]{Lemma} 
\newtheorem{prp}[dfn]{Proposition} 
\newtheorem{thm}[dfn]{Theorem} 
\newtheorem{rmk}[dfn]{Remark} 
\newtheorem{cor}[dfn]{Corollary} 
\newtheorem{ex}[dfn]{Example}

\begin{document}

\bibliographystyle{plain}

\author{Katrin Wehrheim 
\thanks{wehrheim@math.ethz.ch; supported by Swiss National Science Foundation grant 21-64937.01;
} }

\title{Banach space valued Cauchy-Riemann equations with totally real boundary conditions}

\maketitle

\begin{abstract}
The main purpose of this paper is to give a general regularity result for Cauchy-Riemann
equations in complex Banach spaces with totally real boundary conditions.
The usual elliptic $L^p$-regularity results hold true under one crucial assumption:
The totally real submanifold has to be modelled on an $L^p$-space or a closed subspace
thereof.

Secondly, we describe a class of examples of such totally real submanifolds, namely
gauge invariant Lagrangian submanifolds in the space of connections over a Riemann surface.
These pose natural boundary conditions for the anti-self-duality equation on $4$-manifolds
with a boundary space-time splitting, leading towards the definition of a Floer homology for 
$3$-manifolds with boundary, which is the first step in a program by Salamon for the
proof of the Atiyah-Floer conjecture.
The principal part of such a boundary value problem is an example of a Banach space valued 
Cauchy-Riemann equation with totally real boundary condition.
\end{abstract}

\section{Introduction}
A complex Banach space is a Banach space $X$ equipped with a complex structure,
i.e.\ $J\in{\rm End}\,X$ that satisfies $J^2=-\one$.
The Cauchy-Riemann equation for a map $u:\Om\to X$ on a domain $\Om\subset\R^2$
with coordintes $(s,t)$ is $\pd_s u + J\pd_t u = 0$. We will also study the 
equation with an inhomogeneous term on the right hand side.
As in the finite dimensional case, totally real boundary conditions are natural 
for this Cauchy-Riemann equation.
A Banach submanifold $\cL\subset X$ is called totally real with respect to the complex structure 
$J$ if for all $x\in\cL$ one has the direct sum decomposition
$$
X=\rT_x\cL \oplus J\,\rT_x\cL .
$$ 
Let $\Om\subset\H$ be a compact 2-dimensional submanifold in the half space
$$
\H:=\{ (s,t)\in\R^2 \st t\geq 0 \}.
$$
i.e.\ $\Om$ has smooth boundary that might intersect $\pd\H=\{t=0\}$.
We will consider Cauchy-Riemann equations for maps $u:\Om\to X$ that satisfy
totally real boundary conditions on the boundary part $\pd\Om\cap\pd\H$.
Fix an inhomogeneous term $G:\Om\to X$, a family $J:\Om\to{\rm End}\,X$ of
complex structures on $X$, and let $\cL\subset X$ be a Banach submanifold that is
totally real with respect to $J_{s,t}$ for all $(s,t)\in\Om$.
Then we study the following boundary value problem for $u:\Om\to X$,
\begin{equation}\label{holo bvp intro}
\left\{\begin{array}{l}
\pd_s u + J_{s,t}\pd_t u = G ,\\
u(s,0) \in\cL \quad\forall (s,0)\in\pd\Om\cap\pd\H .
\end{array}
\right.
\end{equation}
The Cauchy-Riemann equation itself is linear, but for the linearization of this 
boundary value problem one has to linearize the boundary conditions.
So fix a path $x:\R\to\cL$, then we will also study the Cauchy-Riemann equation 
with linearized totally real boundary conditions for $u:\Om\to X$,
\begin{equation}\label{lin bvp intro}
\left\{\begin{array}{l}
\pd_s u + J_{s,t}\pd_t u = G ,\\
u(s,0) \in \rT_{x(s)}\cL \quad\forall (s,0)\in\pd\Om\cap\pd\H .
\end{array}
\right.
\end{equation}
In this case, there also is a weak formulation of the boundary value problem.
We denote by $X^*$ the dual space of $X$ and denote by $J^*\in{\rm End}\,X^*$ 
the dual operator of the complex structure $J\in{\rm End}\,X$.
Then the weak formulation of (\ref{lin bvp intro}) for $u:\Om\to X$ is
$$
\int_\Om \la u \,,\, \pd_s\psi + \pd_t (J^*\psi)  \ra
\;=\; - \int_\Om \la G \,,\, \psi \ra 
$$
for all $\psi\in\cC^\infty(\Om,X^*)$ with $\supp\psi\subset{\rm int}\,\Om$ and
$\psi(s,0)\in (J(s,0)\rT_{x(s)}\cL)^\perp$ for all $(s,0)\in\pd\Om\cap\pd\H$.
In order to obtain regularity results for any of the above boundary value
problems we make the following crucial assumption.

\begin{description}
\item[$\mathbf{(H_p)}$]
Throughout we suppose that the totally real submanifold $\cL\subset X$ is -- 
as a Banach manifold -- modelled on a closed subspace $Y\subset Z$ of an 
$L^p$-space $Z= L^p(M,\R^m)$ for some $p>1$, $m\in\N$, and a closed manifold $M$.
\end{description}

To show that this assumption still allows $\cL$ to be modelled on a wide variety of 
Banach spaces, we give the following examples.

\begin{ex}  \hspace{1mm} \label{Hp ex}\\
\vspace{-5mm}
\begin{enumerate}
\item 
Every finite dimensional space $\R^m$ is isometric to the subspace of 
constants in $L^p(M,\R^m)$ for ${\rm Vol}\,M=1$.

\item
The Sobolev space $W^{\ell,p}(M)$ (and thus every closed subspace thereof)
is isomorphic to a closed subspace of $L^p(M,\R^m)$.

To see this, choose vector fields $X_1,\ldots,X_k\in\G(\rT M)$ that span $\rT_x M$
for all $x\in M$.
Then the map $u \mapsto (u,\nabla_{X_1}u, \ldots, \nabla_{X_k}^\ell u)$
running through all derivatives of $u$ in the direction of the $X_i$ 
up to order $\ell$ gives an isomorphism between $W^{\ell,p}(M)$ 
and a closed subspace of $L^p(M,\R^m)=:Z$.
\item
Finite products of closed subspaces in $L^p(M_i,\R^{m_i})$ are
isometric to a closed subspace of $L^p(\bigcup M_i,\R^{\max\{m_i\}})$.
\end{enumerate} 
\end{ex}

Our first main theorem gives regularity results and estimates for solutions of 
(\ref{holo bvp intro}) depending on the regularity of the inhomogeneous term in the 
Cauchy-Riemann equation.
Here and throughout the interior of $\Om$ is defined with respect to the topology of 
$\H$, so ${\rm int}\,\Om$ still contains $\pd\Om\cap\pd\H$. 
We use the notation $\N=\{1,2,\ldots\}$. \\

\begin{thm}  \label{holo regularity}
Fix $1<p<\infty$ and a compact subset $K\subset{\rm int}\,\Om$.
Let $\cL\subset X$ be a Banach submanifold that satisfies $(H_p)$.
\begin{enumerate}
\item
Fix $k\in\N$ and let
$$
q:= \left\{\begin{array}{cl}
 p &;\text{if}\;k\geq 2, \\
2p &;\text{if}\;k=1 .
\end{array} \right.
$$
Suppose that $u\in W^{k,q}(\Om,X)$ solves (\ref{holo bvp intro}) for
$G\in W^{k,q}(\Om,X)$ and with a family $J\in W^{k+1,\infty}(\Om,{\rm End}\,X)$ of complex 
structures on $X$, with respect to which $\cL$ is totally real.
Then $u\in W^{k+1,p}(K,X)$.
\item
Let $J_0\in\cC^\infty(\Om,{\rm End}\,X)$ be a smooth family of complex structures 
on $X$, with respect to which $\cL$ is totally real.
Let $u_0\in\cC^\infty(\Om,X)$ be such that $u_0(s,0)\in\cL$ for all 
$(s,0)\in\pd\Om\cap\pd\H$.
Then there exists a constant $\d>0$ with the following significance:
For every constant $c$ and for every $k\in\N$ 
there exists a constant $C$ such that the following holds:
If $u$, $G$, and $J$ satisfy the hypotheses of (i) and
\begin{align*}
\|u-u_0\|_{L^\infty(\Om,X)}\leq\d, \qquad
&\|J-J_0\|_{L^\infty(\Om,{\rm End}\,X)}\leq\d , \\
&\|J-J_0\|_{W^{k+1,\infty}(\Om,{\rm End}\,X)}\leq c ,
\end{align*}
then
$$
\|u-u_0\|_{W^{k+1,p}(K,X)} \leq 
C \bigl( 1+ \|G\|_{W^{k,q}(\Om,X)} + \|u-u_0\|_{W^{k,q}(\Om,X)} \bigr).
$$
\end{enumerate}
\end{thm}

Firstly note the special form of this theorem for $k=1$ (for $k\geq 2$ one has $q=p$).
If $u$ and $G$ are $W^{1,p}$-regular, then one can only deduce $W^{2,\frac p2}$-regularity
of $u$ due to nonlinearities introduced by the coordinates. Moreover, this requires the 
assumption $(H_{\frac p2})$ that $\cL$ is modelled on an $L^{\frac p2}$-space.

Secondly, note that the $u_0$ in (ii) satisfies the Lagrangian boundary condition 
but is not a solution of the Cauchy-Riemann equation. 
It will be required as reference for the construction of coordinates near $\cL$ that 
straighten out the boundary condition but do not depend on the solution $u$ and hence 
allow to deduce an estimate for $u$.
In order that the constant in the estimate becomes independent of the complex structure $J$,
this construction moreover requires that $J$ is $\cC^0$-close to a fixed family $J_0$ 
of complex structures.
The $W^{k+1,\infty}$-bound on the complex structure is only required in order to obtain uniform 
constants.

Moreover, for fixed $k\in\N$ in theorem~\ref{holo regularity} it would actually suffice to have 
$W^{k+1,\infty}$-regularity of $u_0$ and $J_0$. For the interior regularity and estimates
it even suffices to have $W^{k,\infty}$-regularity and bounds on $J$ since one does not need
to write $u$ in coordinates that are adapted to the boundary condition and hence depend on $J$.
This is the same situation as in the finite dimensional case, c.f.\ \cite{MS}.
Finally, the submanifold $\cL$ need only be totally real with respect to $J_{s,t}$ for 
$(s,t)\in\pd\Om\cap\pd\H$. Since this is an open condition, it is then automatically totally 
real in a neighbourhood of $\pd\Om\cap\pd\H$.\\

The second main result concerns the linearized boundary value problem 
(\ref{lin bvp intro}). We use its weak formulation to state the following regularity
result.\\

\begin{thm}   \label{holo regularity 2}
Fix $1<p<\infty$, a compact subset $K\subset{\rm int}\,\Om$, and a Banach submanifold 
$\cL\subset X$ that satisfies $(H_p)$. 
Fix a path $x\in W^{1,\infty}(\R,\cL)$ in $\cL$ and let 
$J\in W^{1,\infty}(\Om,{\rm End}\,X)$ be a family of complex 
structures on $X$, with respect to which $\cL$ is totally real.
Then there is a constant $C$ such that the following holds:

Suppose that $u\in L^p(\Om,X)$ and that there exists a constant $c_u$ such that for 
all $\psi\in W^{1,\infty}(\Om,X^*)$ with $\supp\psi\subset{\rm int}\,\Om$ and
$\psi(s,0)\in (J(s,0)\rT_{x(s)}\cL)^\perp$ for all $(s,0)\in\pd\Om\cap\pd\H$
$$
\left| \int_\Om \la u \,,\, \pd_s\psi + \pd_t (J^*\psi)  \ra \right|
\;\leq\; c_u \|\psi\|_{L^{p^*}(\Om,X^*)} .
$$
Then $u\in W^{1,p}(K,X)$ and
$$
\|u\|_{W^{1,p}(K,X)} \leq C \bigl( c_u + \|u\|_{L^p(\Om,X)} \bigr).
$$
\end{thm}

For strong solutions of the linearized boundary value problem (\ref{lin bvp intro})
the more suitable formulation of theorem~\ref{holo regularity 2} is the following estimate. \\

\begin{cor}  \label{holo regularity 3}
In the setting of theorem~\ref{holo regularity 2} there exists a constant $C$ 
such that the following holds:
Suppose that $u\in W^{1,p}(\Om,X)$ satisfies $u(s,0)\in\rT_{x(s)}\cL$ for all
$(s,0)\in\pd\Om\cap\pd\H$, then 
$$
\|u\|_{W^{1,p}(K,X)} \leq 
C \bigl(\|\pd_s u + J \pd_t u\|_{L^p(\Om,X)} + \|u\|_{L^p(\Om,X)} \bigr).
$$
\end{cor}

A first application of the above results is the elliptic theory for anti-self-dual instantons with 
Lagrangian boundary conditions. It is developed in \cite{W ell}, where theorem~\ref{holo regularity} 
is used to obtain nonlinear regularity and compactness results, whereas theorem~\ref{holo regularity 2} 
enters in the Fredholm theory.
Since the Fredholm theory is performed on a special compact model domain, we make the following remark.

\begin{rmk} \rm 
Theorem~\ref{holo regularity 2} and corollary~\ref{holo regularity 3} 
remain true when $\R$ is replaced by $S^1$, i.e.\ when one considers compact 
domains $K\subset{\rm int}\,\Om\subset S^1\times[0,\infty)$ 
in the half cylinder and a loop $x\in\cC^\infty(S^1,\cL)$.

To see this, identify $S^1\cong\R/\Z$, identify $K$ with a compact subset
$K'\subset\H$ in $[0,1]\times[0,\infty)$, and periodically extend $x$ and $u$ 
for $s\in [-1,2]$. Then $u$ is defined and satisfies the weak equation on some 
open domain $\Om'\subset\H$ such that $K'\subset{\rm int}\,\Om'$, 
so theorem~\ref{holo regularity 2} and corollary~\ref{holo regularity 3} apply. 
These assert regularity and estimates on $K'$ and hence also on $K$.
\end{rmk}

We now proceed to describe a class of examples, to which the above regularity
theory for the Cauchy-Riemann equation can be applied.

A symplectic Banach space $(Z,\o)$ consists of a Banach space $Z$ and a symplectic
structure $\o$, that is a nondegenerate,\footnote{
Nondegeneracy means that for all $z\in Z\setminus\{0\}$ there exists a $y\in Z$ such 
that $\o(z,y)\neq 0$.
}
skewsymmetric, bilinear form $\o:Z\times Z \to\R$. 
In the finite dimensional case there always exists an $\o$-compatible complex 
structure $J\in{\rm End}\,Z$, i.e.\ a complex structure such that 
$\o(\cdot,J\cdot)$ defines a positive definite inner product on $Z$.
In the case of an infinite dimensional Banach space this is not necessarily true.
If an $\o$-compatible complex structure exists, then the norm on $Z$ that is induced by the 
inner product will be bounded but not necessarily complete. The completion of $Z$ with
respect to that norm is then a complex Hilbert space.
In the example below, this Hilbert space will always be the same -- only the complex structure
varies.

Our example of a symplectic Banach space will be the space of connections over a Riemann
surface $\S$. We restrict the discussion to the trivial $\rG$-bundle over $\S$, where
$\rG$ is a compact Lie group.\footnote{
The discussion directly generalizes to nontrivial bundles, where the connections can be
described as $1$-forms with values in an associated bundle.}
Then the space of $L^p$-regular connections is given by the $L^p$-regular $1$-forms with 
values in the Lie algebra $\cg$ of $\rG$. We denote this space
$$
\cA^{0,p}(\S)=L^p(\S;\rT^*\S\otimes\cg) .
$$ 
(For more details on gauge theory and the notation see section~\ref{weakly flat} and 
\cite{W}.)
For $p\geq 2$ the Banach space $\cA^{0,p}(\S)$ is equipped with the symplectic structure
\begin{equation} \label{omega}
\o(\a,\b) = \int_\S \la \a\wedge\b \ra  \qquad\forall \a,\b\in L^p(\S;\rT^*\S\otimes\cg) .
\end{equation}
Moreover, for $p>2$ the gauge group $\cG^{1,p}(\S)=W^{1,p}(\S,\rG)$ acts on $\cA^{0,p}(\S)$ 
by 
$$
u^*A = u^{-1}A u + u^{-1}\rd u \qquad \forall A\in\cA^{0,p}(\S), u\in\cG^{1,p}(\S) .
$$
This gauge action leaves $\o$ invariant. So throughout we assume $p>2$.
Now for any metric on $\S$, the Hodge $*$ operator induces an $\o$-compatible 
complex structure on $\cA^{0,p}(\S)$. The associated inner product is the $L^2$-inner
product of $\cg$-valued $1$-forms, and the completion of $\cA^{0,p}(\S)$ with respect to 
the induced norm is always $L^2(\S,\rT^*\S\otimes\cg)$.

We call a Banach submanifold $\cL\subset (\cA^{0,p}(\S),\o)$ Lagrangian if it is isotropic, 
i.e.\ $\o|_{\cL}\equiv 0$, and if $\rT_A\cL$ is maximal for all $A\in\cL$.
By the latter we mean that for all $\a\in\cA^{0,p}(\S)$ the following implication holds:
$$
\forall\b\in\rT_A\cL \quad \o(\a,\b)=0  \quad\Longrightarrow\quad \a\in\rT_A\cL .
$$ 
In section \ref{weakly flat} we will introduce the space of weakly flat $L^p$-connections
$\cA^{0,p}_{\rm flat}(\Si)$. In particular, we prove that every weakly flat connection is 
gauge equivalent to a smooth connection.
Then we shall show in section~\ref{spaces} that a gauge invariant Lagrangian submanifold of 
$\cA^{0,p}(\S)$ that also satisfies $\cL\subset\cA^{0,p}_{\rm flat}(\Si)$ is automatically 
totally real with respect to the Hodge $*$ operator for any metric on $\S$, 
i.e.\ for all $A\in\cL$
$$
\cA^{0,p}(\S) = \rT_A\cL \oplus * \rT_A\cL .
$$
Moreover, such Lagrangian submanifolds satisfy the assumption $(H_p)$ for theorems~\ref{holo regularity} 
and \ref{holo regularity 2}.
The assumptions of gauge invariance and flatness also ensure that the Lagrangian submanifold
$\cL$ descends to a Lagrangian submanifold in the (singular) symplectic manifold
$M_\S=\cA^{0,p}_{\rm flat}(\S)/\cG^{1,p}(\S)$, the moduli space of gauge equivalence 
classes of flat connections.
The latter can be viewed as symplectic quotient, as was first observed by \cite{AB}.
Note that both $M_\S$ and the quotient $\cL/\cG^{1,p}(\S)$ are allowed to have singularities.
These do not enter the discussion since we will be working in the total space.

Now a pseudoholomorphic curve $u:\Om \to M_\S$ with Lagrangian boundary conditions 
on $\pd\Om\cap\pd\H$ lifts to a solution $B :\Om\times\S \to \rT^*\S\otimes\cg$ 
of the boundary value problem
\begin{equation}\label{pseudohol bvp}
\left\{\begin{array}{rl}
F_B &\hspace{-3mm}= 0 , \\
\pd_s B + *\pd_t B &\hspace{-3mm}= \rd_B \P + * \rd_B \Psi ,\\
B|_{(s,0)\times\S} &\hspace{-3mm}\in\cL \quad\forall (s,0)\in \pd\Om\cap\pd\H .
\end{array}\right.
\end{equation}
Here $\P,\Psi:\Om\times\S \to \cg$ are determined by the solution $B$.
For given $\P,\Psi$, the above boundary value problem without the first equation
is a Cauchy-Riemann equation with totally real boundary conditions as studied in this paper.

Changing the first equation in (\ref{pseudohol bvp}) to 
$*F_B = \pd_t\P -\pd_s\Psi + [\P,\Psi]$ leads to the the anti-self-duality 
equation for the connection $A=\P\ds + \Psi\dt + B$ on $\Om\times\S$
with Lagrangian boundary conditions,
\begin{equation}\label{ASD bvp}
\left\{\begin{array}{l}
*F_A + F_A = 0,\\
A|_{(s,0)\times\S} \in\cL \quad\forall (s,0)\in\pd\Om\cap\pd\H .
\end{array}\right.
\end{equation}
This boundary value problem arises naturally from the Chern-Simons $1$-form on a 
$3$-manifold $Y$ with boundary $\S$ :
This $1$-form becomes closed and it is in fact the differential of the (multivalued) 
Chern-Simons functional, when it is restricted to the space $\cA(Y,\cL)$ of connections $A$
on $Y$ with Lagrangian boundary conditions $A|_\S\in\cL$.
Now the gradient flow lines of the Chern-Simons functional are just the solutions of
(\ref{ASD bvp}) in a special gauge.

It is a program by Salamon \cite{Sa1} to use the boundary value problem (\ref{ASD bvp}) 
to define a Floer homology $\HF^{\rm inst}_*(Y,L)$ for $3$-manifolds $Y$ with boundary $\pd Y=\S$
and Lagrangian submanifolds $L=\cL/\cG^{1,p}(\S)\subset M_\S$, i.e.\ a generalized Morse 
homology for the Chern-Simons functional on $\cA(Y,\cL)$.
As a first indication for the wellposedness of (\ref{ASD bvp}) we prove in 
corollary~\ref{bc approx} that every $W^{1,p}$-regular connection satisfying the boundary 
condition in (\ref{ASD bvp}) can be approximated by smooth connections satisfying the 
same boundary condition.
The elliptic theory for the definition of this Floer homology is set up in \cite{W ell},
where the regularity theorems~\ref{holo regularity} and \ref{holo regularity 2} play a crucial role.
Another approach to the definition of a Floer homology for $3$-manifolds with boundary was
introduced by Fukaya \cite{Fu}. This also uses Lagrangian boundary conditions, but the
construction is restricted to the case of nontrivial bundles, in which case the quotient 
$\cL/\cG^{1,p}(\S)$ is smooth.\\

Finally, a concrete example of a totally real submanifold in a complex Banach space is
given in lemma~\ref{LY}.
Let $\S=\pd Y$ be the boundary of a handle body $Y$ and consider the $L^p$-closure 
of the set of smooth flat connections on $\S$ that can be extended to a flat connection on $Y$,
$$
\cL_Y \,:=\; {\rm cl}\, \bigl\{ A\in\cA_{\rm flat}(\S) \st \exists 
                              \tA\in\cA_{\rm flat}(Y) : \tA|_\S=A \bigr\} 
\;\subset\;\cA^{0,p}(\Si).
$$
This is a Lagrangian submanifold and it is gauge invariant and contained in the space of
flat connections, so as above it also is totally real with respect to the Hodge operator
as complex structure.

These submanifolds occur in the Atiyah-Floer conjecture for homology $3$-spheres as follows:
A Heegard splitting $Y=Y_0\cup_\S Y_1$ of a homology $3$-sphere $Y$ into two handlebodies 
$Y_0$ and $Y_1$ with common boundary $\S$ gives rise to two Lagrangian submanifolds
$\cL_{Y_i}\subset\cA^{0,p}(\Si)$ in the space of $\SU(2)$-connections.
One then has a symplectic Floer homology $\HF^{\rm symp}_*(M_\Si,L_{Y_0},L_{Y_1})$
for the quotients $L_{Y_i}:=\cL_{Y_i}/\cG^{1,p}(\S)\subset M_\S$.
(This is generated by the intersection points of the Lagrangian submanifolds and the 
boundary operator arises from counting pseudoholomorphic strips with Lagrangian 
boundary conditions, i.e.\ solutions of a boundary value problem like (\ref{pseudohol bvp}).)
It was conjectured by Atiyah \cite{A1} and Floer that this should be isomorphic to the instanton
Floer homology $\HF^{\rm inst}_*(Y)$, the generalized Morse homology for the Chern-Simons
functional on the space of $\SU(2)$-connections on $Y$.
Now the program by Salamon \cite{Sa1} is to establish this isomorphism in two steps via the 
intermediate $\HF^{\rm inst}_*([0,1]\times\Si,L_{Y_0}\times L_{Y_1})$ by adiabatic 
limit type arguments similar to \cite{DS}. These adiabatic limits will again require 
elliptic estimates for boundary value problems including a Cauchy-Riemann equation with
totally real boundary conditions as studied in this paper.\\

This paper is organized as follows:
In section \ref{regu} we prove theorems~\ref{holo regularity} and \ref{holo regularity 2} and 
corollary~\ref{holo regularity 3}.
Section \ref{weakly flat} is of preliminary nature: We introduce the notion of a
weakly flat connection, prove the fundamental regularity result for weakly flat 
connections, and discuss the moduli space of flat conections over a Riemann surface.
Section \ref{spaces} deals with gauge invariant Lagrangian submanifolds in the 
space of connections. We establish their basic properties and prove the approximation
result mentioned above. Moreover, we show that the $\cL_Y$ are indeed examples of
Lagrangian and totally real submanifolds.\\

I would like to thank Dietmar Salamon for his constant help and encouragement in pursueing 
this project.

\section{Regularity}
\label{regu}

In this section we prove the regularity theorems~\ref{holo regularity}, \ref{holo regularity 2}, 
and corollary \ref{holo regularity 3}.
Let $\Om\subset\H$ be a compact 2-dimensional submanifold of the half space.
Consider a Banach space $X$ with a family $J:\Om\to{\rm End}\,X$ of complex structures.
Let $\cL\subset X$ be a Banach submanifold that satisfies $(H_p)$, i.e.\ it is modelled on a 
closed subspace $Y\subset Z$ of an $L^p$-space $Z= L^p(M,\R^m)$ for some $p>1$, $m\in\N$, 
and a closed manifold $M$, and suppose that $\cL$ is totally real with respect to
all $J_{s,t}$ for $(s,t)\in\Om$.
Then we consider maps $u:\Om\to X$ that solve the boundary value problem 
(\ref{holo bvp intro}), restated here:
\begin{equation}\label{holo bvp}
\left\{\begin{array}{l}
\pd_s u + J_{s,t}\pd_t u = G ,\\
u(s,0) \in\cL \quad\forall (s,0)\in\pd\Om\cap\pd\H .
\end{array}
\right.
\end{equation}
The idea for the proof of theorem~\ref{holo regularity} is to straighten out the boundary condition 
by going to local coordinates in $Y\times Y$ near $u(s,0)\in X$ such that $Y\times\{0\}$ corresponds 
to the submanifold $\cL$ and the complex structure becomes standard along $Y\times\{0\}$. 
For theorem~\ref{holo regularity 2}, concerning the linearization of (\ref{holo bvp}), one 
chooses $\R$-dependent coordinates for $X$ that identify $Y\times\{0\}$ with $\rT_{x(s)}\cL$ along the 
path $x:\R\to\cL$.
Then the boundary value problem (\ref{holo bvp}) or its linearization yields Dirichlet and 
Neumann boundary conditions for the two components of $u$ and one can use regularity results 
for the Laplace equation with such boundary conditions.

However, there are two difficulties.
Firstly, by straightening out the totally real submanifold, the complex structure 
$J$ becomes explicitly dependent on $u$, so one has to deal carfully with 
nonlinearities in the equation.
Secondly, this approach requires a Cald\'eron-Zygmund inequality for functions 
with values in a Banach space. In general, the Cald\'eron-Zygmund inequality
is only true for values in Hilbert spaces. However, due to the assumption that $\cL$ 
is modelled on an $L^p$-space, we only need the $L^p$-inequality for functions with values 
in $L^p$-spaces.
In that case, the Cald\'eron-Zygmund inequality holds, as can be seen by 
integrating over the real valued inequality.
This will be made precise in the following lemma, in which (i),(iii) are regularity results 
for the homogeneous Dirichlet problem and (ii),(iv) concern the Neumann problem with possibly 
inhomogeneous boundary conditions.
In (i),(ii) the minimum regularity of $u$ is $W^{1,p}$ -- in the case of lower regularity one 
has to use the weak formulation in (iii), (iv).
We abbreviate $\laplace:=\rd^*\rd$ and denote by $\n$ the outer unit normal to $\pd\Om$.
We write $Z^*$ for the dual space of any Banach space $Z$ and write $\la\cdot,\cdot\ra$ for 
the pairing of $Z$ and $Z^*$.
The Sobolev spaces of Banach space valued functions considered below are all
defined as completions of the smooth functions with respect to the respective Sobolev norm.
Moreover, we use the notation 
\begin{align*}
\cC^\infty_\d(\Om,Z^*) 
&:= \{ \psi\in\cC^\infty(\Om,Z^*) \st \psi|_{\pd\Om}=0 \},  \\
\cC^\infty_\n(\Om,Z^*) 
&:= \{ \psi\in\cC^\infty(\Om,Z^*) \st \tfrac{\pd\psi}{\pd\n}\bigr|_{\pd\Om}=0 \}.
\end{align*}

\pagebreak

\begin{lem} \label{Banach space regularity}
Fix $1<p<\infty$ and $k\in\N$ and 
let $\Om$ be a compact Riemannian manifold with boundary.
Let $Z=L^p(M)$ for some closed manifold $M$.
Then there exists a constant $C$ such that the following holds.
\begin{enumerate}
\item Let $f\in W^{k-1,p}(\Om,Z)$ and suppose that $u\in W^{k,p}(\Om,Z)$
solves
$$
\int_\Om \la u \,,\, \laplace\psi \ra 
= \int_\Om \la f \,,\, \psi \ra 
\qquad\forall \psi\in\cC^\infty_\d(\Om,Z^*) .
$$
Then $u\in W^{k+1,p}(\Om,Z)$ and
$\; \|u\|_{W^{k+1,p}} \leq C \|f\|_{W^{k-1,p}} $.
\item Let $f\in W^{k-1,p}(\Om,Z)$, $g\in W^{k,p}(\Om,Z)$, 
and suppose that $u\in W^{k,p}(\Om,Z)$ solves
$$
\int_\Om \la u \,,\, \laplace\psi \ra 
= \int_\Om \la f \,,\, \psi \ra + \int_{\pd\Om} \la g \,,\, \psi \ra 
\qquad\forall \psi\in\cC^\infty_\n(\Om,Z^*) .
$$
Then $u\in W^{k+1,p}(\Om,Z)$ and
$$
\|u\|_{W^{k+1,p}} \leq C \bigl( \|f\|_{W^{k-1,p}} + \|g\|_{W^{k,p}} 
                              + \|u\|_{L^p} \bigr).
$$
\item Suppose that $u\in L^p(\Om,Z)$ and there exists a constant $c_u$ such that 
$$
\left| \int_{\Om\times M} u \cdot \laplace_\Om \psi \,\right|
\;\leq\; c_u \|\psi\|_{W^{1,p^*}(\Om,Z^*)} 
\qquad\forall\psi\in\cC^\infty_\d(\Om\times M) .
$$
Then $u\in W^{1,p}(\Om,Z)$ and $\;\|u\|_{W^{1,p}(\Om,Z)} \leq C c_u $.
\item Suppose that $u\in L^p(\Om,Z)$ and there exists a constant $c_u$ such that
$$
\left| \int_{\Om\times M} u \cdot \laplace_\Om \psi \,\right|
\;\leq\; c_u \|\psi\|_{W^{1,p^*}(\Om,Z^*)}
\qquad\forall\psi\in\cC^\infty_\n(\Om\times M) .
$$
Then $u\in W^{1,p}(\Om,Z)$ and
$$
\|u\|_{W^{1,p}(\Om,Z)} \leq C \bigl( c_u + \|u\|_{L^p(\Om,Z)}\bigr).
$$
If moreover $\int_\Om u = 0$ then in fact 
$\quad\|u\|_{W^{1,p}(\Om,Z)} \leq C c_u $.
\end{enumerate}
\end{lem}

The key to the proof of (i) and (ii) is the fact that the functions $f$ and $g$ 
can be approximated not only by smooth functions with values in the Banach
space $L^p(M)$, but by smooth functions on $\Om\times M$.

\begin{lem} \label{approx}
Let $\Om$ be a compact manifold (possibly with boundary), let $M$ be a closed
manifold, let $1<p,q<\infty$, and $k,\ell\in\N_0$. Then the following holds.

\begin{enumerate}
\item 
$\cC^\infty(\Om\times M)$ is dense in $W^{k,q}(\Om,W^{\ell,p}(M))$.
\item 
A function $u\in W^{k,q}(\Om,W^{\ell,p}(M))$ with zero boundary values $u|_{\pd\Om}=0$ 
can be approximated by $u^\n\in\cC^\infty(\Om\times M)$ with $u^\n|_{\pd\Om\times M}=0$.
\item
If $\ell p>\dim M$ and $z\in M$, then a function $u\in W^{k,q}(\Om,W^{\ell,p}(M))$ with 
$u(\cdot,z)=0 \in W^{k,q}(\Om)$ can be approximated by $u^\n\in\cC^\infty(\Om\times M)$ with 
$u^\n(\cdot,z)\equiv 0$.
\end{enumerate}
\end{lem}

\noindent
{\bf Proof of lemma \ref{approx}: } \\
We first prove (i).
By definition $\cC^\infty(\Om,W^{\ell,p}(M))$ is dense in 
$W^{k,q}(\Om,W^{\ell,p}(M))$. So we fix 
$g\in\cC^\infty(\Om,W^{\ell,p}(M))$ and show that in every 
$W^{k,q}(\Om,W^{\ell,p}(M))$-neighbourhood of $g$ there exists a
$\tilde g\in\cC^\infty(\Om\times M)$.
Firstly, we prove this in the case $k=0$ for closed manifolds $M$ as well as
in the following case (that will be needed for the proof in the case $k\geq 1$): 
$M=\R^n$, $g$ is supported in $\Om\times V$ and $\tilde g$ is required to
have support in $\Om\times U$ for some open bounded domains $V,U\subset\R^n$
such that $\overline V\subset U$.

Fix $\d>0$. Since $\Om$ is compact one finds a finite covering 
$\Om=\bigcup_{i=1}^N U_i$ by neighbourhoods $U_i$ of $x_i\in\Om$ such that
$$
\| g(x)-g(x_i) \|_{W^{\ell,p}(M)} \leq \tfrac\d 2
\qquad\forall x\in U_i .
$$
Next, choose $g_i\in\cC^\infty(M)$ such that 
$\| g_i - g(x_i) \|_{W^{\ell,p}(M)} \leq \tfrac\d 2$.
In the case $M=\R^n$ one has $\supp g(x_i)\subset V$ and hence can
choose $g_i$ such that it is supported in $U$ (e.g.\ using mollifiers with 
compact support).
Then choose a partition of unity $\sum_{i=1}^N \p_i = 1$ by
$\p_i\in\cC^\infty(\Om,[0,1])$ with $\supp\p_i\subset U_i$.
Now one can define $\tilde g \in\cC^\infty(\Om\times M)$ by
$$
\tilde g(x,z) := \sum_{i=1}^N \p_i(x) g_i(z)
\qquad\forall x\in\Om, z\in M .
$$
In the case $M=\R^n$ this satisfies $\supp\tilde g \subset \Om\times U$ as
required. Moreover,
\begin{align*}
\|\tilde g - g\|_{L^q(\Om,W^{\ell,p}(M))}^q
&= \int_\Om \bigl\| \tsum_{i=1}^N \p_i(g_i - g) \bigr\|_{W^{\ell,p}(M)}^q \\
&\leq \int_\Om \bigl( \tsum_{i=1}^N \p_i \cdot
       \sup_{x\in U_i} \| g_i - g(x) \|_{W^{\ell,p}(M)} \bigr)^q \\
&\leq \int_\Om \d^q
\;=\; \d^q \, {\rm Vol}\,\Om .
\end{align*}
Thus we have proven the lemma in the case $k=0$. For $k\geq 1$ this method does 
not work since one picks up derivatives of the cutoff functions $\p_i$.
Instead, one has to use mollifiers and the result for $k=0$ on $M=\R^n$.

So we assume $k\geq 1$, fix $g\in\cC^\infty(\Om,W^{\ell,p}(M))$ and pick some $\d>0$.
Let $M=\bigcup_{i=1}^N \P_i(U_i)$ be an atlas with bounded open domains 
$U_i\subset\R^n$ and charts $\P_i: U_i \to M$.
Let $V_i\subset \overline V_i \subset U_i$ be open sets such that still
$M=\bigcup_{i=1}^N \P_i(V_i)$. Then there exists a partition of unity
$\sum_{i=1}^N \psi_i\comp\P_i^{-1}=1$ by 
$\psi_i\in\cC^\infty(\R^n,[0,1])$ such that $\supp\psi_i\subset V_i$.
Now $g=\sum_{i=1}^N g_i{\scriptscriptstyle \circ}({\rm id}_\Om \times \P_i^{-1})$ with
$$
g_i(x,y) = \psi_i(y) \cdot g(x,\P_i(y)) \qquad\forall x\in\Om, y\in U_i .
$$
Here $g_i\in\cC^\infty(\Om,W^{\ell,p}(\R^n))$ is extended by $0$ outside of
$\supp g_i \subset \Om \times V_i$, and it suffices to prove that each of
these functions can be approximated in $W^{k,q}(\Om,W^{\ell,p}(\R^n))$ by 
$\tilde g_i \in \cC^\infty(\Om\times \R^n)$ with 
$\supp\tilde g_i \subset \Om \times U_i$.
So drop the subscript $i$ and consider $g\in\cC^\infty(\Om,W^{\ell,p}(\R^n))$ 
that is supported in $\Om\times V$, where $V,U\subset\R^n$ are open bounded 
domains such that $\overline V \subset U$.

Let $\s_\ep(y)=\ep^{-n}\s(y/\ep)$ be a family of compactly supported mollifiers 
for $\ep>0$, i.e.\ $\s\in\cC^\infty(\R^n,[0,\infty))$ such that 
$\supp\s\subset B_1(0)$ and $\int\s = 1$. Then for all $\ep>0$ define 
$\tilde g_\ep\in\cC^\infty(\Om\times\R^n)$ by
$$
\tilde g_\ep (x,y) := [\s_\ep * g(x,\cdot)](y)
\qquad\forall x\in\Om, y\in\R^n .
$$
Firstly, $\supp\s_\ep\subset B_\ep(0)$, so for sufficiently small $\ep>0$
the support of $\tilde g_\ep$ lies within $\Om\times U$.
Secondly, we abbreviate for $j\leq k$, $m\leq\ell$
$$
f_{j,m}\,:=\;\nabla_\Om^j\nabla_{\R^n}^m g 
\;\in\;\cC^\infty(\Om,L^p(\R^n)) ,
$$
which are supported in $\Om\times V$. Then
\begin{align*}
\|\tilde g_\ep - g\|_{W^{k,q}(\Om,W^{\ell,p}(\R^n))}^q
&= \sum_{j\leq k} \int_\Om \bigl\| \nabla_\Om^j\bigl(\s_\ep*g(x,\cdot) 
                            - g(x,\cdot) \bigr) \bigr\|_{W^{\ell,p}(\R^n)}^q \\
&\leq (\ell+1)^{\frac qp} \sum_{j\leq k} \sum_{m\leq\ell} \int_\Om 
\bigl\| \s_\ep * f_{j,m}(x,\cdot) - f_{j,m}(x,\cdot) \bigr\|_{L^p(\R^n)}^q .
\end{align*}
Now use the result for $k=0$ on $M=\R^n$ (with values in a vector bundle) 
to find $\tilde f_{j,m}\in\cC^\infty(\Om\times\R^n)$ supported
in $\Om\times U$ such that
$$
\|\tilde f_{j,m} - f_{j,m}\|_{L^q(\Om,L^p(\R^n))}\leq\d .
$$
Then for all $x\in\Om$ and sufficiently small $\ep>0$ the functions
$\s_\ep * \tilde f_{j,m}(x,\cdot)$ are supported in some fixed bounded 
domain $U'\subset\R^n$ containing $U$.
Moreover, the $\tilde f_{j,m}$ are Lipschitz continuous, hence one finds
a constant $C$ (depending on the $\tilde f_{j,m}$, i.e.\ on $g$ and $\d$) 
such that for all $x\in\Om$
\begin{align*}
& \bigl\| \s_\ep * \tilde f_{j,m}(x,\cdot) 
                 - \tilde f_{j,m}(x,\cdot) \bigr\|_{L^p(\R^n)}^p  \\
&= \int_{U'} \Bigl| \int_{\R^n}  \s_\ep(y'-y) \bigl( \tilde f_{j,m}(x,y') 
                     - \tilde f_{j,m}(x,y) \bigr) \rd^n y'\Bigr|^p \rd^n y \\
&\leq \int_{U'} \Bigl( \int_{\R^n}  \s_\ep(y'-y) 
      \sup_{|y-y'|\leq\ep} |\tilde f_{j,m}(x,y') - \tilde f_{j,m}(x,y)|
      \,\rd^n y' \Bigr)^p \rd^n y \\
&\leq   {\rm Vol}\, U'(C\ep)^p.
\end{align*}
Now use the fact that the convolution with $\s_\ep$ is continuous with respect to the 
$L^p$-norm, $\|\s_\ep * f \|_p\leq \|f\|_p$ (see e.g.\ \cite[Lemma 2.18]{Adams}) to estimate

\begin{align*}
&\int_\Om 
 \bigl\| \s_\ep * f_{j,m}(x,\cdot) - f_{j,m}(x,\cdot) \bigr\|_{L^p(\R^n)}^q \\
&\leq \int_\Om \Bigl(
 \bigl\| \s_\ep * \bigl(       f_{j,m}(x,\cdot) 
                      - \tilde f_{j,m}(x,\cdot) \bigr) \bigr\|_{L^p(\R^n)} 
+\bigl\| f_{j,m}(x,\cdot) - \tilde f_{j,m}(x,\cdot) \bigr\|_{L^p(\R^n)} \\
&\qquad\quad +\bigl\| \s_\ep * \tilde f_{j,m}(x,\cdot) 
                - \tilde f_{j,m}(x,\cdot) \bigr\|_{L^p(\R^n)} \Bigr)^q \\
&\leq 2\cdot 3^q \bigl\| f_{j,m} - \tilde f_{j,m} \bigr\|_{L^q(\Om,L^p(\R^n))}^q  
      + 3^q \,{\rm Vol}\,\Om \,({\rm Vol}\,U)^{\frac qp} (C\ep)^q 
\quad\leq\; 3\cdot 3^q \d^q .
\end{align*}
Here we have chosen 
$0<\ep\leq C^{-1}({\rm Vol}\,\Om)^{-\frac 1q}({\rm Vol}\,U)^{-\frac 1p}\,\d$.
Thus we obtain
$$
\|\tilde g_\ep - g\|_{W^{k,q}(\Om,W^{\ell,p}(\R^n))} \leq 
3(\ell+1)^{\frac 1p}(3(k+1)(\ell+1))^{\frac 1q}\, \d .
$$
This proves (i).
To show (ii) one first approximates in $\cC^\infty(\Om,W^{\ell,p}(M))$ with zero
boundary values and then mollifies on $M$ as in (i) as follows.

In case $k=0$ the boundary condition is meaningless, but the approximation with zero
boundary values can be done elementary by cutting off in small neighbourhoods of the boundary.
For $k\geq 1$ consider a local chart of $\Om$ in $[0,1]\times\R^n$ such that $\{t=0\}$ 
corresponds to the boundary, where $t$ denotes the $[0,1]$-coordinate. 
Let $f\in W^{k,q}([0,1]\times\R^n,Z)$ for any vector space $Z$ with $f|_{t=0}=0$ and 
compact support.
Let $\s_\ep$ be mollifiers on $\R^n$ as above, then
$f_\ep(t,\cdot):=\s_\ep * f(t,\cdot)$ defines $f_\ep\in \cC^\infty(\R^n,W^{k,q}([0,1],Z))$
for all $\ep>0$. One checks that $\|f_\ep- f\|_{W^{k,q}([0,1]\times\R^n)}\to 0$ as $\ep\to 0$.
We choose the $\s_\ep$ with compact support, then the $f_\ep$ are also compactly
supported and hence have finite $W^{k,q}([0,1],W^{\ell,q}(\R^n))$-norm for any $\ell\in\N$.
Moreover, note that still $f_\ep|_{t=0}=0$.
In order to approximate $f_\ep$ with zero boundary values one chooses $\ell=k$, then (i)
gives a smooth approximation $g^\n\to\pd_t f_\ep$ in the
$W^{k-1,q}([0,1],W^{k,q}(\R^n))$-norm.
Now $f_\ep^\n(t,x):=\int_0^t g^\n(\t,x)\rd\t$ defines functions in 
$\cC^\infty([0,1]\times\R^n,Z)$ that vanish at $t=0$ and approximate $f_\ep$ in the
$W^{k,q}([0,1],W^{k,q}(\R^n))$-norm, which is even stronger than the 
$W^{k,q}$-norm on $[0,1]\times\R^n$.


Finally, to prove (iii), we choose an approximation by $u^\n\in\cC^\infty(\Om\times M)$.
Then $u^\n(\cdot,z)\to 0$ in $W^{k,q}(\Om)$ since the evaluation at $z$ is a continuous map 
$W^{\ell,p}(M)\to\R$. 
Now $u^\n - u^\n(\cdot,z)\in\cC^\infty(\Om\times M)$ still converges to $u$ in 
$W^{k,q}(\Om,W^{\ell,p}(M))$ but it vanishes at $z$.
\QED

In the case $q=p$ lemma~\ref{approx} provides the continuous inclusion 
$$
W^{k,p}(\Om,W^{\ell,p}(M)) \subset W^{\ell,p}(M,W^{k,p}(\Om))
$$ 
since the norms on these spaces are identical.\footnote{
The spaces are actually equal. The proof requires an extension of the 
approximation argument to manifolds with boundary. We do not carry this
out here because we will only need this one inclusion.
}
Moreover, for $p=q$ and $k=\ell=0$ the lemma identifies $L^p(\Om,L^p(M)) = L^p(\Om\times\Si)$ 
as the completion of $\cC^\infty(\Om\times M)$ under the $L^p$-norm.\\

\noindent
{\bf Proof of lemma \ref{Banach space regularity} (i) and (ii) : } \\
We first give the proof of the regularity for the inhomogeneous Neumann problem
(ii) in full detail; (i) is proven in complete analogy
-- using the regularity theory for the Laplace equation on $\R$-valued 
functions with Dirichlet boundary condition instead of the Neumann condition.

Fix $f\in W^{k-1,p}(\Om,Z)$, $g\in W^{k,p}(\Om,Z)$, and let 
$f^i,g^i\in\cC^\infty(\Om\times M)$ be approximating sequences given by
lemma~\ref{approx}. 
Testing the weak equation with $\psi\equiv\a$ for all $\a\in Z^*$
implies $\int_\Om f + \int_{\pd \Om} g = 0$ and thus 
$h^i :=\int_\Om f^i + \int_{\pd\Om} g^i \to 0$ in $Z$ as $i\to\infty$, so 
one can replace the $f^i$ by $f^i-h^i / {\rm Vol}\,\Om \in\cC^\infty(M,Z)$ to achieve
$$
\int_\Om f^i(\cdot,y) + \int_{\pd\Om} g^i(\cdot,y) = 0 
\qquad \forall y\in M, i\in\N .
$$
Now for each $y\in M$ there exist unique solutions 
$u^i(\cdot,y)\in\cC^\infty(\Om)$ of 
\[
\left\{\begin{array}{rl}
\laplace u^i (\cdot,y) &\hspace{-2.5mm}= f^i(\cdot,y) , \\
\tfrac{\pd}{\pd\n} u^i(\cdot,y)\bigr|_{\pd\Om} 
&\hspace{-2.5mm}= g^i(\cdot,y)\bigr|_{\pd\Om}, \\
\tint_\Om u^i(\cdot,y) &\hspace{-2.5mm}= 0  .
\end{array}\right.
\]
For each of these Laplace equations with Neumann boundary conditions one
obtains an $L^p$-estimate for the solution, see e.g.\ \cite[Proposition~A.1]{W ell}
and \cite[Theorem 3.1]{W} for the existence. The constant can be chosen
independently of $y\in M$ since it varies continuously with $y$ and $M$ is compact. 
Then integration of those estimates yields (with different constants $C$)
\begin{align*}
\| u^i \|_{W^{k+1,p}(\Om,Z)}^p
&= \int_M \bigl\| u^i \bigr\|_{W^{k+1,p}(\Om)}^p \\
&\leq \int_M C \bigl( \| f^i \|_{W^{k-1,p}(\Om)} 
                    + \| g^i \|_{W^{k,p}(\Om)} \bigr)^p \\
&\leq C \bigl( \| f^i \|_{W^{k-1,p}(\Om,Z)}
             + \| g^i \|_{W^{k,p}(\Om,Z)}    \bigr)^p .
\end{align*}
Here one uses the crucial fact that
$L^p(\Om,L^p(M))\subset L^p(M,L^p(\Om))$ with identical norms.
(Note that this is not the case if the integrability indices over $\Om$ and $M$
are different.)
Similarly, one obtains for all $i,j\in\N$
$$
\| u^i - u^j \|_{W^{k+1,p}(\Om,Z)}
\leq C \bigl( \| f^i - f^j \|_{W^{k-1,p}(\Om,Z)}
            + \| g^i - g^j \|_{W^{k,p}(\Om,Z)}    \bigr) .
$$
So $u^i$ is a Cauchy sequence and hence converges to some
$\tilde u\in W^{k+1,p}(\Om,Z)$. 
Now suppose that $u\in W^{k,p}(\Om,Z)$ solves the weak Neumann equation for
$f$ and $g$, then we claim that in fact $u=\tilde u + c \in W^{k+1,p}(\Om,Z)$,
where $c\in Z$ is given by
$$
c(y) := \frac 1{{\rm Vol}\,\Om}
        \int_\Om \bigl( u(\cdot,y) - \tilde u(\cdot,y)  \bigr)
\qquad\forall y\in M .
$$
In order to see that indeed $c\in L^p(M)=Z$ and that for some constant $C$ one has  
$\|c\|_{L^p(M)}\leq C (\|u\|_{L^p(\Om,Z)} + \|\tilde u\|_{L^p(\Om,Z)} )$
note that lemma~\ref{approx} yields the continuous inclusion 
$W^{k,p}(\Om,L^p(M))\subset L^p(M,W^{k,p}(\Om))\subset L^p(M,L^1(\Om))$.
To establish the identity  $u=\tilde u + c$, we first note that for all 
$\p\in\cC^\infty(M)\subset Z^*$
$$
\int_\Om \la \tilde u + c - u \,,\, \p \ra 
\;=\; \int_M \p \cdot \Bigl( c - \int_\Om (\tilde u - u) \Bigr) 
\;=\; 0 .
$$
Next, for any $\p\in\cC^\infty(\Om\times M)$ let 
$$
\p_0:=\frac 1{{\rm Vol}\,\Om}\int_\Om\p \quad\in\;\cC^\infty(M) .
$$
Then one finds $\psi\in\cC^\infty_\n(\Om\times M)$
such that $\p=\laplace_\Om\psi + \p_0$.
(There exist unique solutions $\psi(\cdot,y)$ of the Neumann problem 
for $\p(\cdot,y)-\p_0(y)$, and these depend smoothly on $y\in M$.)
So we find that for all $\p\in\cC^\infty(\Om\times M)$, abbreviating $\laplace_\Om=\laplace$
\begin{align*}
\int_\Om \la u - \tilde u - c \,,\, \p \ra 
&= \int_\Om \la u \,,\, \laplace\psi \ra 
    - \int_\Om \la \tilde u + c \,,\, \laplace\psi \ra 
    + \int_\Om \la u - \tilde u - c \,,\, \p_0 \ra     \\
&= \lim_{i\to\infty}\left( \int_\Om \la f - \laplace u^i \,,\, \psi \ra 
         + \int_{\pd\Om} \la g - \tfrac{\pd u^i}{\pd\n} \,,\, \psi \ra \right)
\;=\; 0 .
\end{align*}
This proves $u=\tilde u + c\in W^{k+1,p}(\Om,Z)$ and the 
estimate for $u^i$ yields in the limit
\begin{align*}
\|u\|_{W^{k+1,p}(\Om,Z)} 
&\leq \|\tilde u\|_{W^{k+1,p}(\Om,Z)} 
    + ({\rm Vol}\,\Om)^{\frac 1p} \| c \|_{L^p(M)} \\
&\leq C \bigl( \|f\|_{W^{k-1,p}(\Om,Z)} + \|g\|_{W^{k,p}(\Om,Z)} 
             + \|u\|_{L^p(\Om,Z)} \bigr).
\end{align*}
This finishes the proof of (ii), and analogously of (i).
\QED \\

\noindent
{\bf Proof of lemma \ref{Banach space regularity} (iii) and (iv) : } \\
Let $u\in L^p(\Om,Z)$ be as supposed in (iii) or (iv), where $Z=L^p(M)$ and thus 
$Z^*=L^{p^*}(M)$. 
Then we have $u\in L^p(\Om\times M)$ and the task is to prove that
$\rd_\Om u$ also is of class $L^p$ on $\Om\times M$.
So we have to consider $\int_{\Om\times M} u\cdot \rd_\Om^*\t$
for $\t\in\cC^\infty_\d(\Om\times M,\rT^*\Om)$
(which are dense in $L^{p^*}(\Om\times M,\rT^*\Om)$).
In the case (iii) one finds for any such smooth family $\t$ of $1$-forms on $\Om$
a smooth function $\psi\in\cC^\infty_\d(\Om\times M)$ such that 
$\rd_\Om^*\t=\laplace_\Om\psi$.
Then there is a constant $C$ such that for all $y\in M$
(see e.g. \cite[Proposition~A.1]{W ell} 
$$
\|\psi(\cdot,y)\|_{W^{1,p^*}}
\;\leq\; C\,\|\laplace_\Om\psi(\cdot,y)\|_{(W^{1,p})^*}
\;\leq\; C\,\|\t(\cdot,y)\|_{p^*} .
$$ 
In the case (iv) one similarly finds $\psi\in\cC^\infty_\n(\Om\times M)$ such that 
$\rd_\Om^*\t=\laplace_\Om\psi$ and 
$\|\psi(\cdot,y)\|_{W^{1,p^*}}\leq C\,\|\t(\cdot,y)\|_{p^*}$ 
for all $y\in M$ and some constant $C$.
(Note that $\int_\Om \rd_\Om^*\t \equiv 0$ since $\t$ vanishes on $\pd\Om\times M$
and we have used e.g.\ \cite[Theorems 2.2,2.3']{W}.)
In both cases we can thus estimate for all $\t\in\cC^\infty_\d(\Om\times M,\rT^*\Om)$
using the assumption
\begin{align*}
\left| \int_{\Om\times M}  u \cdot \rd_\Om^*\t  \,\right|
&\;=\; \left| \int_{\Om\times M}  u \cdot \laplace_\Om\psi \,\right| 
\;\leq\; c_u \left( \int_M \|\psi\|_{W^{1,p^*}(\Om)}^{p^*} \right)^{\frac 1{p^*}}\\
&\;\leq\; C c_u \left( \int_M \|\t\|_{L^{p^*}(\Om)}^{p^*} \right)^{\frac 1{p^*}} 
\;\leq\; C c_u \| \t \|_{L^{p^*}(\Om\times M)} .
\end{align*}
Now in both cases the Riesz representation theorem (e.g.\ \cite[Theorem 2.33]{Adams}) 
asserts that $\int_{\Om\times M} u \cdot \rd_\Om^*\t = \int_{\Om\times M} f \cdot \t $
for all $\t$ with some $f\in L^p(\Om\times M)$.
This proves the $L^p$-regularity of $\rd_\Om u$ and yields the estimate 
$$
\|\rd_\Om u\|_{L^p(\Om\times M)}\leq C c_u .
$$
In the case (iii), one can moreover deduce $u|_{\pd\Om}=0$. Indeed, partial integration
in the weak equation gives for all $\psi\in\cC^\infty_\d(\Om\times M)$
\begin{align*}
\left| \int_{\pd\Om\times M} u \cdot \tfrac{\pd\psi}{\pd\n} \right|
&= \left| \int_{\Om\times M}  u \cdot \laplace_\Om\psi \;\;
        - \int_{\Om\times M} \la \rd_\Om u \,,\, \rd_\Om\psi \ra   \right| \\
&\leq (c_u + \|\rd_\Om u\|_{L^p(\Om\times M)} ) \|\psi\|_{W^{1,p^*}(\Om\times M)}.
\end{align*}
For any given $g\in\cC^\infty(\pd\Om\times M)$ one now finds 
$\psi\in\cC^\infty(\Om\times M)$ with $\psi|_{\pd\Om}=0$ and 
$\frac{\pd\psi}{\pd\n}=g$, and these can be chosen such that 
$\|\psi\|_{W^{1,p^*}}$ becomes arbitrarily small. Then one obtains
$\int_{\pd\Om\times M} u \, g = 0$ and thus $u|_{\pd\Om}=0$.
Thus in the case~(iii) one finds a constant $C'$ such that
$$
\|u\|_{W^{1,p}(\Om,Z)}^p 
\;=\; \int_M \| u \|_{W^{1,p}(\Om)}^p  
\;\leq\; C' \int_M \| \rd_\Om u \|_{L^p(\Om)}^p  
\;=\; C' \| \rd_\Om u \|_{L^p(\Om\times M)}^p  ,
$$
which finishes the proof of (iii).

In case (iv) with the additional assumption $\int_\Om u = 0$ one also has a constant 
$C'$ such that 
$ \| u(\cdot,y) \|_{W^{1,p}(\Om)}\leq C' \| \rd_\Om u (\cdot,y)\|_{L^p(\Om)} $ 
for all $y\in M$ and thus
$$
\|u\|_{W^{1,p}(\Om,Z)}
\;\leq\; C' \| \rd_\Om u \|_{L^p(\Om\times M)}  
\;\leq\; C'Cc_u .
$$
In the general case (iv) one similarly has
$$
\|u\|_{W^{1,p}(\Om,Z)}^p 
\;=\; \| \rd_\Om u \|_{L^p(\Om\times M)}^p + \| u \|_{L^p(\Om\times M)}^p  
\;\leq\; \bigl( C c_u  + \| u \|_{L^p(\Om,Z)} \bigr) ^p   .
$$
\QED

The proof of theorem~\ref{holo regularity} will moreover use the following 
quantitative version of the implicit function theorem. This is proven e.g.\ 
in \cite[Proposition A.3.4]{MS} by a Newton-Picard method.
(Here we only need the special case $x_0=x_1=0$.)

\begin{prp} \label{qIFT}
Let $X$ and $Y$ be Banach spaces and let $U\subset Y$ be a neighbourhood of $0$.
Suppose that $f:\cU\to X$ is a continuously differentiable map such that 
$\rd_0 f : Y\to X$ is bijective.
Then choose constants $c\geq\|(\rd_0 f)^{-1}\|$ and $\d>0$ such that
$B_\d(0)\subset\cU$ and 
$$
\| \rd_y f - \rd_0 f \| \leq \tfrac 1 {2c} \qquad\forall y \in B_\d(0)  .
$$
Now if $\|f(0)\|\leq\frac\d{4c}$ then there exists a unique solution 
$y\in B_\d(0)$ of $f(y)=0$.
Moreover, this solution satisfies 
$$
\|y\|\leq 2 c \|f(0)\|.
$$
\end{prp}

\noindent
{\bf Proof of theorem~\ref{holo regularity} : } \\
Let $z_0\in\cL$ and let $J_0\in{\rm End}\,X$ be a complex structure with respect to which
$\cL$ is totally real.
Choose a Banach manifold chart $\p:V\to\cL$ from a neighbourhood $V\subset Y$ 
of $0$ to a neighbourhood of $\p(0)=z_0$.
Then one obtains a Banach submanifold chart of $\cL\subset X$ from a neighbourhood 
$\cW \subset Y\times Y$ of zero to a ball $B_\ep(z_0)\subset X$ around $z_0$,
$$
\Th:
\begin{array}{ccc}
\cW & \overset{\sim}{\longrightarrow} & B_\ep(z_0)  \\
(v_1,v_2) & \longmapsto & \p(v_1) + J_0 \rd_{v_1}\p (v_2) .
\end{array}
$$
To see that this is indeed a diffeomorphism for sufficiently small $\cW$ and 
$\ep>0$ we just check that $D:=\rd_{(0,0)}\Th = \rd_0\p \oplus J_0 \rd_0\p$ is an 
isomorphism.
This is since $\rd_0\p:Y\to\rT_{z_0}\cL$ is an isomorphism and so is the map
$\rT_{z_0}\cL \times \rT_{z_0}\cL \to X$ given by the splitting
$X=\rT_{z_0}\cL \oplus J_0\,\rT_{z_0}\cL$.
The size $\ep>0$ of the chart can be quantified by proposition~\ref{qIFT} as follows.
For the maps $f=\Th - x : V\times V \to X$ one finds constants $c=\|D^{-1}\|$ and $\d>0$ 
independently of $x\in X$ such that $B_\d(0)\subset V\times V$ and
$$
\| \rd_y f - \rd_0 f \| \;=\; \| \rd_y \Th - \rd_0 \Th \| \leq \tfrac 1 {2c} 
\qquad\forall y \in B_\d(0)  .
$$
Then for $\|x-z_0\|\leq \frac\d{4c}=:\ep$ one obtains a unique $y=\Th^{-1}(x)$
in $B_\d(0)$.

Next, if one replaces $z_0$ and $J_0$ by $z\in\cL$ and a complex structure $J\in{\rm End}\,X$ 
in sufficiently small neighbourhoods of $z_0$ and $J_0$ respectively, then one still obtains 
a Banach submanifold chart $\Th:\cW\to B_\ep(z)$ with $\Th(0)=z$.
Here $\cW$ varies with $(z,J)$, but one can choose a uniform $\ep>0$. 
This is since one can find uniform constants $c$ and $\d$ in proposition~\ref{qIFT}.
(The map $\Th$ varies with $z$ via the chart map $\p_z:V\to\cL$, 
$v\mapsto \p(\p^{-1}(z)+v)$ that is defined for sufficiently small $V$ and that
satisfies $\p_z(0)=z$.)
Moreover, one obtains the uniform estimate
\begin{equation}\label{Th est}
\|\Th^{-1}(x)\|_{Z\times Z} \leq C \|x-z\|_X \qquad \forall x\in B_\ep(z) .
\end{equation}
(Recall that $Y$ is a closed subspace of the Banach space $Z$, so the norm on
$Y$ is induced by the norm on $Z$.)

Now consider a solution $u\in W^{k,q}(\Om,X)$ of (\ref{holo bvp intro}) for some
$G\in W^{k,q}(\Om,X)$ and $J\in W^{k+1,\infty}(\Om,{\rm End}\,X)$ as in 
theorem~\ref{holo regularity}~(i).
Fix any $(s_0,0)\in K$ and let $z\equiv u(s_0,0)\in\cL$.
Then the above construction of the coordinates $\Th$ can be done for all $J=J_{s,t}$ 
with $(s,t)\in U$ for a neighbourhood $U\subset\Om$ of $(s_0,0)$.
Thus one obtains a $W^{k+1,\infty}$-family of chart maps for $(s,t)\in\overline{U}$, 
$$ 
\Th_{s,t}: Y\times Y \supset \cW_{s,t} \overset{\sim}{\longrightarrow} 
           B_\ep(z_{s,t}) .
$$
Recall that $u$ is either of class $W^{1,2p}$ or of class $W^{k,p}$ with $k\geq 2$ and 
$p>1$. On the $2$-dimensional domain $\Om$, the Sobolev embeddings thus ensure that $u$ is
continous. So on a possibly even smaller neighbourhood $U$ of $(s_0,0)$ the map $u$ 
can be expressed in local coordinates,
$$
u(s,t) = \Th_{s,t}(v(s,t)) \qquad\forall (s,t)\in U ,
$$
where $v\in W^{k,q}(U,Z\times Z)$. This follows from the fact that the composition of the
$W^{k+1,\infty}$-map $\Th^{-1}$ with a $W^{k,q}$-map $u$ is again $W^{k,q}$-regular if $kq>2$ 
(see e.g.\ \cite[Lemma B.8]{W}).
Moreover, $v$ actually takes values in $\cW\subset Y\times Y$. 

In order to obtain the estimate in (ii), the map $\Th$ has to be constructed independently 
of $u$ and $J$, using the fixed $u_0$ and $J_0$. In that case let 
\hbox{$z_{s,t}:=u_0(s,0)$}, which 
is welldefined on a small neighbourhood $U$ of $(s_0,0)\in K\subset{\rm int}\,\Om$. Then 
the coordinates $\Th_{s,t}$ are defined for all $(s,t)\in U$ and for all complex structures 
in a sufficiently small neighbourhood of $J_0(s_0,0)$.
In particular, $\Th_{s,t}$ is defined for all $J=J_{s,t}$ with $(s,t)\in U$, provided that 
$J\in W^{k+1,\infty}(\Om,{\rm End}\,X)$ satisfies the assumption $\|J-J_0\|_{L^\infty}\leq\d$. 
Here one again makes sufficiently small choices of $U$ and $\d>0$.
Thus one obtains a $W^{k+1,\infty}$-family of chart maps $\Th_{s,t}$ as above that now also
satisfy the uniform estimate (\ref{Th est}) for all $(s,t)\in U$, where the constant $C$
only depends on $u_0$ and $J_0$.
Now in order to again express $u$ in local coordinates, choose $U$ even smaller such that 
$u_0(s,t)\in B_{\frac \ep 2}(x_0)$ for all $(s,t)\in U$ and let $\d\leq\frac \ep 2$.
Then every $u\in W^{k,q}(\Om,X)$ that satisfies $\|u-u_0\|_{L^\infty(\Om,X)}\leq\d$ 
can be written $u = \Th\comp v$ as above.
Now integration of (\ref{Th est}) together with the fact that all derivatives
of $\Th^{-1}$ up to order $k$ are bounded (due to the $W^{k,\infty}$-bound on $J$) 
yields the estimate 
$$
\|v\|_{W^{k,q}(U,Z\times Z)}
\leq C \|u-u_0\|_{W^{k,q}(U,X)} .
$$
Here and in the following $C$ denotes any constant that is independent of the 
specific choices of $J$ and $u$ in the fixed neighbourhoods of $J_0$ and $u_0$, 
however, it may depend on $c$ and $k$.

In the coordinates constructed above, the boundary value problem (\ref{holo bvp}) now becomes
\begin{equation} \label{v bvp}
\left\{\begin{array}{l}
\pd_s v + I \pd_t v = f , \\
v_2(s,0)=0 \quad\forall s\in\R .
\end{array}
\right.
\end{equation}
with $v=(v_1,v_2)$ and
\begin{align*}
f &= (\rd_v\Th)^{-1} \bigl( G - \pd_s\Th (v) - J \pd_t\Th(v) \bigr)
  \;\in\; W^{k,q}(U,Y\times Y) , \\
I &= (\rd_v\Th)^{-1} J \rd_v\Th  
\qquad\qquad\qquad\qquad\;\;\in\; W^{k,q}(U,{\rm End}(Y\times Y)) .
\end{align*}
Note the following difficulty: The complex structure $I$ now explicitly depends on the solution
$v$ of the equation (\ref{v bvp}) and thus is only $W^{k,q}$-regular. This cannot be avoided
when straightening out the Lagrangian boundary condition.
However, one obtains one more simplification of the boundary value problem:
$\Th$ was constructed such that one obtains the standard complex structure along $\cL$. 
Indeed, for all $(s,0)\in U$ using that $J^2=-\one$
\begin{align*}
I(s,0) 
&= (\rd_{(v_1,0)}\Th)^{-1} J \rd_{(v_1,0)}\Th       
= \bigl(\rd_{v_1}\p \oplus J \rd_{v_1}\p \bigr)^{-1}  
 J \bigl(\rd_{v_1}\p \oplus J \rd_{v_1}\p \bigr)      \\ 
&=
\left(\begin{array}{rr}
       0 & -\one \\
       \one & 0
\end{array}\right) 
\;=:\; I_0.
\end{align*}
Moreover, in case (ii) one has the following estimates on $U$:
\begin{align*}
\| I \|_{W^{k,q}} &\leq C, \\
\| f \|_{W^{k,q}} &\leq C \bigl( \|G\|_{W^{k,q}} + \|u-u_0\|_{W^{k,q}} \bigr) .
\end{align*}
So for every boundary point $(s_0,0)\in K\cap\pd\H$ we have rewritten the 
boundary value problem (\ref{holo bvp}) over some neighbourhood $U\subset\Om$.
Now for the compact set $K\subset\Om$ one finds a covering
$K\subset V \cup \bigcup_{i=1}^N U_i $ by finitely many such neighbourhoods 
$U_i$ at the boundary and a compact domain $V\subset\Om\setminus\pd\Om$ away 
from the boundary.
Note that the $U_i$ can be replaced by interior domains $\tU_i$ (that intersect
$\pd U_i$ only on $\pd\H$) that together with $V$ still cover $K$.
We will establish the regularity and estimate for $u$ on all domains
$\tU_i$ near the boundary and on the remaining domain $V$ separately.
So firstly consider a domain $U_i$ near the boundary and drop the subscript 
$i$. After possibly replacing $U$ by a slightly smaller domain one can assume 
that $U$ is a manifold with smooth boundary and still $\tU\cap\pd U \subset \pd\H$.
The task is now to prove the regularity and estimate for $u=\Th\comp v$ on $\tU$ from 
$(\ref{v bvp})$. 

Since $\Th_{s,t} : Y\times Y \to X$ are smooth maps in $W^{k+1,\infty}$-dependence on 
\hbox{$(s,t)\in U$}, it suffices to prove that $v\in W^{k+1,p}(\tU,Z\times Z)$ 
with the according estimate. (One already knows that $v$ takes values
-- almost everywhere -- in $Y\times Y$, so one automatically also obtains
$v\in W^{k+1,p}(\tU,Y\times Y)$.)
For that purpose fix a cutoff function $h\in\cC^\infty(\H,[0,1])$
with $h\equiv 1$ on $\tU$ and $h\equiv 0$ on $\H\setminus U$.
Moreover, this function can be chosen such that $\pd_t h|_{t=0}=0$. 
Note that $h\equiv 0$ on $\pd U \setminus \pd\H$, so $h v_2$ satisfies 
the Dirichlet boundary condition on $\pd U$. 
Indeed, we will see that $hv_2\in W^{k,p}(U,Z)$ solves a weak Dirichlet 
problem.

In the following, $Y^*$ denotes the dual space of $Y$ and we write $\la\cdot,\cdot\ra$ 
for both the pairings between $Y$ and $Y^*$ and between $Y\times Y$ and $Y^*\times Y^*$.
We obtain
$\laplace= - (\pd_s + \pd_t I^*)(\pd_s - I^*\pd_t) + (\pd_t I^*)\pd_s - (\pd_s I^*)\pd_t$
for $I^*\in W^{k,q}(\Om,{\rm End}(Y^*\times Y^*))$
the pointwise dual operator of $I$.
Thus for all $\p\in\cC^\infty(\Om,Y^*\times Y^*)$
\begin{align*}
h \laplace \p
&= - (\pd_s + \pd_t I^*)(\pd_s - I^*\pd_t)(h\p)
  - (\laplace h) \p + 2 (\pd_s h) \pd_s\p + 2 (\pd_t h) \pd_t\p \\
&\quad + (\pd_t I^*)\pd_s(h\p) - (\pd_s I^*)\pd_t(h\p) .
\end{align*}
Hence (\ref{v bvp}) and partial integration (for smooth approximations of $v$, $f$, $I$)
yields
\begin{align}
&\int_U \la h v \,,\, \laplace \p \ra \nonumber\\
&= \int_U \la \pd_s v + I \pd_t v \,,\, (\pd_s - I^*\pd_t)(h\p) \ra 
\nonumber\\
&\quad
 - \int_U \la  (\laplace h) v + 2(\pd_s h) \pd_s v + 2 (\pd_t h) \pd_t v
                + h (\pd_t I) \pd_s v - h (\pd_s I) \pd_t v  \,,\, \p \ra 
\nonumber\\
&\quad
+ \int_{\pd U\cap\pd\H}  \la I v \,,\, (\pd_s - I^*\pd_t)(h\p) \ra
                          + \la h (\pd_s I) v - 2 (\pd_t h) v \,,\, \p \ra \nonumber\\
&= \int_U \la h ( - \pd_s f + I \pd_t f + (\pd_t I) f 
                    - (\pd_t I) \pd_s v  + (\pd_s I) \pd_t v )  \nonumber\\
&\qquad\qquad\qquad\qquad\quad
  - (\laplace h) v - 2(\pd_s h) \pd_s v - 2 (\pd_t h) \pd_t v  \,,\, \p \ra 
\nonumber\\
&\quad
+ \int_{\pd U\cap\pd\H}  \la h \cdot I f \,,\, \p \ra 
+ \int_{\pd U\cap\pd\H}  \la   v \,,\, \pd_t (h\p) \ra
                       + \la I v \,,\, \pd_s (h\p) \ra  \nonumber \\
&= \int_U \la F \,,\, \p \ra + \int_{\pd U}  \la H \,,\, \p \ra 
+ \int_{\pd U\cap\pd\H}  \la v_1 \,,\, \pd_t (h\p_1) + \pd_s (h\p_2) \ra .
\label{weak eqn}
\end{align}
This uses the notation $\p=(\p_1,\p_2)$, the boundary condition 
$v_2|_{t=0}=0$, and the fact that $I|_{t=0}\equiv I_0$.
One then reads off \hbox{$F=(F_1,F_2)\in W^{k-1,p}(U,Y\times Y)$},
$H=(H_1,H_2)\in W^{k,p}(U,Y\times Y)$, and that in case (ii) for some constants $C$
\begin{align*}
\|F\|_{W^{k-1,p}} + \|H\|_{W^{k,p}} 
&\leq C \bigl( \|f\|_{W^{k,q}} + \|I\|_{W^{k,q}}\|f\|_{W^{k,q}} 
            + \|I\|_{W^{k,q}}\|v\|_{W^{k,q}} \bigr) \\
&\leq C \bigl( \|G\|_{W^{k,q}} + \|u-u_0\|_{W^{k,q}} \bigr) .
\end{align*}
We point out that the crucial terms here are $(\pd_s I) \pd_t v$ and $(\pd_t I) \pd_s v$.
In the case $k\geq 3$ the estimate holds with $q=p$ due to the Sobolev embedding
$W^{k-1,p}\cdot W^{k-1,p}\hookrightarrow W^{k-1,p}$. In the case $k=1$ one only has 
$L^{2p}\cdot L^{2p} \hookrightarrow L^p$ and hence one needs $q=2p$ in the above estimate.
In the case $k=2$ the Sobolev embedding $W^{1,p}\cdot W^{1,p}\hookrightarrow W^{1,p}$
requires $p>2$. Let us make this assumption for the moment and finish the proof of this
theorem under the additional assumption $p>2$ in case $k=2$. Then after that we will show
how an iteration of the theorem for the case $k=1$ will give the required $W^{1,p}$-regularity
of $(\pd_s I) \pd_t v$ and $(\pd_t I) \pd_s v$ also in case $1<p\leq 2$.
(Note that this would follow from $W^{2,\tp}$-regularity of $u$ for any $\tp>2$.)

Now in order to obtain a weak Laplace equation for $v_2$ we test the weak equation
(\ref{weak eqn}) with $\p=(\p_1,\p_2)=(0,\pi\comp\psi)$ for $\psi\in\cC^\infty_\d(U,Z^*)$ and 
where $\pi:Z^*\to Y^*$ is the canonical embedding. In that case, both boundary terms 
vanish and one obtains for all $\psi\in\cC^\infty_\d(U,Z^*)$
$$
\int_U \la h v_2 \,,\, \laplace \psi \ra
= \int_U \la F_2 \,,\, \psi \ra .
$$
By lemma~\ref{Banach space regularity}~(i) this weak equation for $hv_2\in W^{k,p}(U,Z)$ 
now implies that $hv_2\in W^{k+1,p}(U,Z)$ and thus $v_2\in W^{k+1,p}(\tU,Z)$.
Moreover, one obtains the estimate
\begin{align*}
\|v_2\|_{W^{k+1,p}(\tU,Z)}
\leq \|hv_2\|_{W^{k+1,p}(U,Z)}  
&\leq C \|F_2\|_{W^{k-1,p}(U,Z)} \\
&\leq C \bigl( \|G\|_{W^{k,q}(\Om,X)} + \|u-u_0\|_{W^{k,q}(\Om,X)} \bigr) .
\end{align*}
To obtain a weak Laplace equation for $v_1$ we test the weak equation (\ref{weak eqn})
with $\p=(\p_1,\p_2)=(\pi\comp\psi,0)$, where $\psi\in\cC^\infty(U,Z^*)$ such that 
\hbox{$\pd_t\psi|_{t=0}=0$}. This makes the second boundary term vanish, so we
obtain for all $\psi\in\cC^\infty_\n(U,Z^*)$
$$
\int_\Om \la h v_1 \,,\, \laplace \psi \ra 
= \int_\Om \la F_1 \,,\, \psi \ra + \int_{\pd\Om} \la H_1 \,,\, \psi \ra .
$$
So we have established a weak Laplace equation with Neumann boundary condition 
for $hv_1$. Now lemma~\ref{Banach space regularity}~(ii) implies that 
$hv_1\in W^{k+1,p}(U,Z)$, hence $v_1\in W^{k+1,p}(\tU,Z)$.
Moreover, one obtains the estimate
\begin{align*}
\|v_1\|_{W^{k+1,p}(\tU,Z)}
&\leq \|hv\|_{W^{k+1,p}(U,Z)}  \\
&\leq C \bigl( \|F_1\|_{W^{k-1,p}(U,Z)} + \|H_1\|_{W^{k,p}(U,Z)}
             + \|h v_1\|_{W^{k,p}(U,Z)} \bigr) \\
&\leq C \bigl( \|G\|_{W^{k,q}(\Om,X)} + \|u-u_0\|_{W^{k,q}(\Om,X)} \bigr) .
\end{align*}
This now provides the regularity and the estimate for $u=\Th\comp v$ on $\tU$ as follows.
We have established that $v : \tU \to Z\times Z$ is a $W^{k+1,p}$-map
that takes values in $\cW\subset Y\times Y$.
All derivatives of $\Th:\Om\times \cW \to X$ up to order $k+1$ are uniformly bounded on $\Om$. 
Hence $u\in W^{k+1,p}(\tU,X)$ and 
\begin{align*}
\|u-u_0\|_{W^{k+1,p}(\tU,X)} 
&\leq C \bigl( 1+ \|v\|_{W^{k+1,p}(\tU,X)} \bigr) \\
&\leq 
C \bigl( 1+ \|G\|_{W^{k,q}(\Om,X)} + \|u-u_0\|_{W^{k,q}(\Om,X)} \bigr).
\end{align*}
For the regularity of $u$ on the domain $V\subset\Om\setminus\pd\Om$ away from 
the boundary one does not need any special coordinates.
As for $U$, one replaces $\Om$ by a possibly smaller domain with smooth boundary.
Moreover, one chooses a cutoff function $h\in\cC^\infty(\H,[0,1])$ such that 
$h|_V\equiv 1$ and that vanishes outside of $\Om\subset\H$ and in a 
neighbourhood of $\pd\Om$.
Then in the same way as for (\ref{weak eqn}) one obtains a weak Dirichlet 
equation. For all $\p\in\cC^\infty_\d(\Om,X^*)$
\begin{align*}
\int_\Om \la h u \,,\, \laplace \p \ra   
&= \int_\Om \la h \bigl( - \pd_s G + J \pd_t G + (\pd_t J) G 
                    - (\pd_t J) \pd_s u  + (\pd_s J) \pd_t u \bigr)  \\
&\qquad\qquad\qquad\qquad\qquad
  - (\laplace h) u - 2(\pd_s h) \pd_s u - 2 (\pd_t h) \pd_t u  \,,\, \p \ra .
\end{align*}
Note that $X\cong Y\times Y\subset Z\times Z$ also is bounded isomorphic to
a closed subspace of an $L^p$-space.
So by lemma~\ref{Banach space regularity} this weak equation implies 
that $hu\in W^{k+1,p}(\Om,X)$, and thus $u\in W^{k+1,p}(V,X)$ with the estimate
$$
\|u\|_{W^{k+1,p}(V,X)}
\leq C \bigl( 1 + \|G\|_{W^{k,q}(\Om,X)} + \|u-u_0\|_{W^{k,q}(\Om,X)} \bigr) .
$$
(Note that here it suffices to have a $W^{k,\infty}$-bound on $J$.)
Thus we have proven the regularity and estimates of $u$ on all parts of the
finite covering $K\subset V\cup \bigcup_{i=1}^N U_i $.
This finishes the proof of the theorem under the additional assumption $p>2$ in case $k=2$.

Finally, let $u\in W^{2,p}(\Om,X)$ and $G\in W^{2,p}(\Om,X)$ be as assumed for 
\hbox{$1<p\leq 2$}. 
Then the task is to establish $W^{2,\tp}$-regularity and -estimates for some $\tp>2$.
This follows from the following iteration.

One starts with $W^{2,p_0}(\Om_0)$-regularity and -estimates for $p_0=p\in (1,2)$ on $\Om_0=\Om$. 
(In case $p=2$ one chooses a smaller value for $p$.)
Now as long as $p_i\leq\frac{4p}{2+p}<2$ and $\Om_i\subset\H$ is compact one has the Sobolev 
embeddings $W^{2,p}(\Om)\hookrightarrow W^{1,\frac{p_i}{2-p_i}}(\Om_i)$ and
$W^{2,p_i}(\Om_i)\hookrightarrow L^{\frac{2p_i}{2-p_i}}(\Om_i)$.
So choose a compact submanifold $\Om_{i+1}\subset\H$ such that $K\subset{\rm int}\,\Om_{i+1}$ 
and $\Om_{i+1}\subset{\rm int}\,\Om_i$. Then the theorem in case $k=1$ with $p$ replaced by
$\frac{p_i}{2-p_i}$ gives regularity and estimates in $W^{2,p_{i+1}}(\Om_{i+1})$ for 
$p_{i+1}=\frac{p_i}{2-p_i}$.
One sees that the sequence $(p_i)$ grows at a rate of at least $\frac 1{2-p}>1$ until it
reaches $p_N\geq \frac{4p}{2+p}$ after finitely many steps.
A further step of the iteration with $p_N=\frac{4p}{2+p}$ then gives 
$W^{2,p_{N+1}}$-regularity and -estimates with $p_{N+1}=\frac{4p}{2-p}>2$ on $\Om_{N+1}:=K$.
This is exactly the regularity for $u$ that still was to be established.
\QED \\

\noindent
{\bf Proof of theorem \ref{holo regularity 2} : } \\
The Banach manifold charts near the path $x:\R\to\cL$ give rise to a $W^{1,\infty}$-path 
of isomorphisms $\p_s:Y\overset{\sim}{\to}\rT_{x(s)}\cL$ for all $s\in\R$.
Together with the family of complex structures 
$J\in W^{1,\infty}(\overline\Om,{\rm End}\,X)$ these give rise to a family
$\Th\in W^{1,\infty}(\overline\Om,{\rm Hom}(Y\times Y,X))$ of bounded isomorphisms
$$
\Th_{s,t}:
\begin{array}{ccc}
Y\times Y & \overset{\sim}{\longrightarrow} & X \\
(z_1,z_2) & \longmapsto & \p_s(z_1) + J_{s,t} \p_s(z_2) .
\end{array}
$$
The inverses of the dual operators of $\Th_{s,t}$ give a family of bounded isomorphisms
$\Th'\in W^{1,\infty}(\overline\Om,{\rm Hom}(Y^*\times Y^*,X^*))$ , 
$$
\Th'_{s,t}:=(\Th_{s,t}^*)^{-1} : \; 
Y^*\times Y^* \; \overset{\sim}{\longrightarrow} \; X^* . 
$$
One checks that for all $(s,t)\in\overline\Om$
$$
\Th_{s,t}^{-1}J_{s,t}\Th_{s,t} 
\;=\; \textstyle {0 \; -\smone \choose \smone \,\;\;\; 0}\;=:\; I_0
\;\in\;{\rm End}(Y\times Y) .
$$
Next, after possibly replacing $\Om$ by a slightly smaller domain that still 
contains $K$ in its interior, one can assume that $\Om$ is a manifold with smooth
boundary.
Then fix a cutoff function $h\in\cC^\infty(\H,[0,1])$ such that $h|_K\equiv 1$
and $\supp h\subset\Om$, i.e.\ $h\equiv 0$ near $\pd\Om\setminus\pd\H$. 
Now let $u\in L^p(\Om,X)$ be given as in the theorem
and express it in the above coordinates as
$u = \Th\comp v$, where $v\in L^p(\Om,Y\times Y)$.
We will show that $v$ satisfies a weak Laplace equation.
For all $\p\in\cC^\infty(\Om,Y^*\times Y^*)$ we introduce
$\psi:=\Th'((\pd_s + I_0\pd_t)\p )\in W^{1,\infty}(\overline{\Om},X^*)$ and calculate
\begin{align*}
\pd_s (h\psi) + \pd_t(J^* h\psi) 
&= h \Th'((\pd_s + I_0^*\pd_t)(\pd_s + I_0\pd_t)\p ) \\
&\quad + (\pd_s h)\psi + \pd_t ( h J^*)\psi 
 + h( \pd_s\Th' + J^*\pd_t\Th' ){\Th'}^{-1}(\psi)  .
\end{align*}
If $\psi(s,0)\in (J(s,0)\rT_{x(s)}\cL)^\perp$ for all 
$(s,0)\in\pd\Om\cap\pd\H$, then $h\psi$ is an admissible test function in the
given weak estimate for $u$ in the theorem and we obtain, denoting all constants by $C$
and using $\Th^*\Th'={\rm id}$,
\begin{align*}
\left| \int_\Om \la h v \,,\, \laplace \p \ra \right| 
&= \left| \int_\Om \la \Th(v) \,,\, 
        h\Th'\bigl((-\pd_s + I_0\pd_t)(\pd_s + I_0\pd_t) \p \bigr) \ra \right|\\
&= \left| \int_\Om \la u \,,\, \pd_s (h\psi) + \pd_t(J^* h\psi)  \ra \right| \\
&\quad + \left|\int_\Om \la u \,,\, (\pd_s h)\psi + \pd_t (h J^*)\psi 
              + h( \pd_s\Th' + J^*\pd_t\Th'){\Th'}^{-1}(\psi) \ra \right|\\
&\leq \bigl( c_u + C \|u\|_{L^p(\Om,X)} \bigr) \|\psi\|_{L^{p^*}(\Om,X^*)} \\
&\leq C \bigl( c_u + \|u\|_{L^p(\Om,X)} \bigr) 
        \|\p\|_{W^{1,p^*}(\Om,Y^*\times Y^*)} .
\end{align*}
Here we used the fact that $J^*$ and $\Th'$ as well as their first derivatives and
inverses are bounded linear operators between $Y^*\times Y^*$ and $X^*$.
This inequality then holds for all $\p=(\p_1,\p_2)$ with 
$\p_1\in\cC^\infty_\n(\Om,Y^*)$ and $\p_2\in\cC^\infty_\d(\Om,Y^*)$ since in
that case $\psi$ is admissible. Indeed, 
 $\psi|_{t=0}=\Th'(\pd_s\p_1-\pd_t\p_2,0)\in (J\rT_x\cL)^\perp$
due to $\Th'(Y^*\times\{0\})=\Th'(I_0(Y\times\{0\}))^\perp=(J\rT_x\cL)^\perp$.

Recall that $Y\subset Z$ is a closed subset of the Banach space $Z$ with the 
induced norm. So one has $v\in L^p(\Om,Z\times Z)$. Let $\pi:Z^*\to Y^*$ 
be the natural embedding, then above inequality holds with 
$\p=(\pi\comp\psi_1,\pi\comp\psi_2)$ for all $\psi_1\in\cC^\infty_\n(\Om,Z^*)$ 
and $\psi_2\in\cC^\infty_\d(\Om,Z^*)$.
Since $\|\pi\comp\psi_i\|_{Y^*} \leq \|\psi_i\|_{Z^*}$ one then obtains for
all such $\Psi=(\psi_1,\psi_2)\in\cC^\infty(\Om,Z^*\times Z^*)$
$$
\left| \int_\Om \la h v \,,\, \laplace \Psi \ra \right| 
\leq C \bigl( c_u + \|u\|_{L^p(\Om,X)} \bigr) 
       \|\Psi\|_{W^{1,p^*}(\Om,Z^*\times Z^*)} .
$$
Now lemma~\ref{Banach space regularity}~(iii) and (iv) asserts the $W^{1,p}$-regularity
of $hv$ and hence one obtains $v\in W^{1,p}(\Om,Z\times Z)$ with the estimate
$$
\|v\|_{W^{1,p}(K,Z\times Z)}
\;\leq\; \|h v\|_{W^{1,p}(\Om,Z\times Z)}
\;\leq\; C\bigl( c_u + \| u \|_{L^p(\Om,X)} + \| v \|_{L^p(\Om,Z\times Z)} \bigr).
$$
For the first factor of $Z\times Z$, this follows from 
lemma~\ref{Banach space regularity}~(iv), in the second factor one uses (iii).
Since it was already known that $v$ takes values in $Y\times Y$ (almost 
everywhere), one in fact has $v\in W^{1,p}(\Om,Y\times Y)$ with the same estimate
as above.
Finally, recall that $u=\Th\comp v$ and use the fact that all derivatives up to first order 
of $\Th$ and $\Th^{-1}$ are bounded to obtain $u\in W^{1,p}(K,X)$ with the claimed 
estimate (using again \cite[Lemma B.8]{W})
$$
\|u\|_{W^{1,p}(K,X)}
\;\leq\; C \|v\|_{W^{1,p}(K,Z\times Z)}  
\;\leq\; C \bigl( c_u + \| u \|_{L^p(\Om,X)} \bigr) .
$$
\QED

\noindent
{\bf Proof of corollary \ref{holo regularity 3}: } \\
Let $u\in W^{1,p}(\Om,X)$ and $\psi\in W^{1,\infty}(\Om,X^*)$ such that
$\supp\psi\subset{\rm int}\,\Om$ and with the boundary conditions
$u(s,0)\in\rT_{x(s)}\cL$ and $\psi(s,0)\in ( J(s,0)\rT_{x(s)}\cL )^\perp$ 
for all $(s,0)\in\pd\Om$. Then one obtains the weak estimate, where the boundary term vanishes,
\begin{align*}
\left| \int_\Om \la u \,,\, \pd_s\psi + \pd_t (J^*\psi) \ra  \right|
&=\left| \int_\Om \la \pd_s u + J \pd_t u \,,\, \psi \ra  
       - \int_{\pd\Om\cap\pd\H} \la J u \,,\, \psi \ra  \right|  \\
&\leq \| \pd_s u + J \pd_t u \|_{L^p(\Om,X)}  \|\psi\|_{L^{p^*}(\Om,X^*)} .
\end{align*}
This holds for all $\psi$ as above, so the estimate follows from theorem~\ref{holo regularity 2}.
\QED

\section{Weakly flat connections}
\label{weakly flat}

In this section we consider the trivial $\rG$-bundle over a closed manifold $\S$ 
of dimension $n\geq 2$. Here $\rG$ is a compact Lie group with Lie algebra $\cg$. 
We recall that $\cg$ is equipped with a Lie bracket $[\cdot,\cdot]$ and a 
$\rG$-invariant inner product $\la\cdot,\cdot\ra$ that moreover satisfy the relation
$$
\la [\x,\e], \z\ra = \la \x, [\e,\z] \ra 
\qquad\forall \x,\z,\e\in\cg .
$$
A (smooth) connection on this bundle is a $\cg$-valued $1$-form $A\in\Om^1(\S;\cg)$.
The exterior derivative $\rd_A$ associated to it is given by 
$\rd_A\e = \rd\e + [A\wedge\e]$ for all $\cg$-valued differential forms $\e$.
Here the Lie bracket indicates how the values of the differential forms are paired.
Now $\rd_A\comp\rd_A$ vanishes if and only if the connection is flat, that is its
curvature $F_A=\rd A + \half [A\wedge A]$ vanishes.
 
Now fix $p>n$ and consider the space $\cA^{0,p}(\Si)=L^p(\S,\rT^*\S\otimes\cg)$
of $L^p$-regular connections. Their curvature is not welldefined, but the flatness
condition can also be formulated weakly:
A connection $A\in\cA^{0,p}(\Si)$ is called {\bf weakly flat} if
\begin{equation}\label{flat}
\int_\Si \la A \,,\, \rd^*\o -\half (-1)^n * [A\wedge *\o] \ra = 0
\qquad\forall\o\in\Om^2(\Si;\cg).
\end{equation}
For sufficiently regular connections one sees by partial integration that (\ref{flat}) 
is equivalent to the connection being flat.
We denote the space of weakly flat $L^p$-connections over $\Si$ by
$$
\cA^{0,p}_{\rm flat}(\Si) 
:= \bigl\{ A\in\cA^{0,p}(\Si) \st A \;\text{satisfies}\; (\ref{flat}) \bigr\}.
$$
One can check that this space is invariant under the action of the gauge
group $\cG^{1,p}(\Si)=W^{1,p}(\S,\rG)$,
$$
u^*A = u^{-1}A u + u^{-1}\rd u \qquad \forall A\in\cA^{0,p}(\S), u\in\cG^{1,p}(\S) .
$$
Note that (\ref{flat}) is welldefined for $p\geq 2$, but $\cG^{1,p}(\Si)$ and its 
action on $\cA^{0,p}(\S)$ are only welldefined for $p>n$, see e.g.\ 
\cite[Appendix B]{W}. 
The next theorem shows that the quotient $\cA^{0,p}_{\rm flat}(\Si)/\cG^{1,p}(\Si)$
can be identified with the usual moduli space of flat connections
$\cA_{\rm flat}(\Si)/\cG(\Si)$ -- smooth flat connections modulo smooth gauge 
transformations.

\begin{thm}  \label{flat regularity} 
For every weakly flat connection $A\in\cA_{\rm flat}^{0,p}(\Si)$ there exists 
a gauge transformation $u\in\cG^{1,p}(\Si)$ such that $u^*A\in\cA_{\rm flat}(\Si)$
is smooth. 
\end{thm}

The proof will be based on the following $L^p$-version of the local slice 
theorem, a proof of which can be found in \cite[Theorem 8.3]{W}.

\begin{prp} \label{weak Coulomb gauge thm} 
Fix a reference connection $\hat A\in\cA^{0,p}(\Si)$.
Then there exists a constant $\d>0$ such that for every $A\in\cA^{0,p}(\Si)$ 
with $\|A-\hat A\|_p \leq \d$ there exists a gauge transformation 
$u\in\cG^{1,p}(\Si)$ such that
\begin{equation}\label{weak Coulomb gauge}
\int_\Si \bigl\langle u^*A-\hat A \,,\, \rd_{\hat A}\e \bigr\rangle = 0 
\qquad\forall \e\in\cC^\infty(\Si,\cg) .
\end{equation}
Equivalently, one has for $v=u^{-1}\in\cG^{1,p}(\Si)$
$$
\int_\Si \bigl\langle v^*\hat A- A \,,\, \rd_A\e \bigr\rangle = 0 
\qquad\forall \e\in\cC^\infty(\Si,\cg) .
$$
\end{prp}

The weak flatness together with the weak Coulomb gauge condition (\ref{weak Coulomb gauge})
form an elliptic system, so theorem~\ref{flat regularity} is then a consequence of the 
regularity theory for the Laplace operator, or the Hodge decomposition of $L^p$-regular 
$1$-forms. \\

\noindent
{\bf Proof of theorem \ref{flat regularity} : } \\
Consider a weakly flat connection $A\in\cA^{0,p}_{\rm flat}(\Si)$.
Let $\d>0$ be the constant from proposition~\ref{weak Coulomb gauge thm} 
for the reference connection $A$ and choose a smooth connection 
$\tA\in\cA(\Si)$ such that $\|\tA - A\|_p \leq \d$.
Then by proposition~\ref{weak Coulomb gauge thm} there exists a gauge transformation 
$u\in\cG^{1,p}(\Si)$ such that
$$
\int_\Si \bigl\langle u^*A- \tA \,,\, \rd_\tA\e \bigr\rangle = 0 
\qquad\forall \e\in\cC^\infty(\Si,\cg) .
$$
Now \cite[Lemma~A.2]{W ell} asserts that $\a:=u^*A-\tA\in L^p(\Si,\rT^*\Si\otimes\cg)$ is in 
fact smooth. 
(By the definition of Sobolev spaces via coordinate charts it suffices to prove 
the regularity and estimate for $\a(X)$, where $X\in\G(\rT \Si)$ is any smooth 
vector field on $\Si$. Alternatively to this lemma -- a consequence of the regularity
theory for the Laplace operator -- one can also deduce the regularity of $\a$ directly
from the regularity of the Hodge decomposition.)
This is due to the weak equations
\begin{align*}
\int_\Si \la \a \,,\, \rd\e \ra 
&= - \int_\Si \la *[\a\wedge *\tA] \,,\, \e \ra
\;\,\qquad\qquad\qquad\forall\e\in\cC^\infty(\Si,\cg), \\
\int_\Si \la \a \,,\, \rd^*\o \ra 
&= - \int_\Si \la \rd\tA + \half [u^*A\wedge u^*A] \,,\, \o \ra 
\qquad\forall\o\in\Om^2(\Si;\cg).
\end{align*}
Firstly, the inhomogeneous terms are of class $L^{\frac p2}$, hence the lemma
asserts $W^{1,\frac p2}$-regularity of $\a$ and $u^*A$.
Now if $p\leq 2n$, then the Sobolev embedding gives $L^{p_1}$-regularity of 
$u^*A$ with $p_1:=\frac{np}{2n-p}$ (in case $p=2n$ one can choose any $p_1>2n$). 
This is iterated to obtain $L^{p_j}$-regularity for the sequence 
$p_{j+1}=\frac{n p_j}{2n-p_j}$ (or any $p_{j+1}>2n$ in case $p_j\geq 2n$) with 
$p_0=p$.
One checks that $p_{j+1}\geq\th p_j$ with $\th=\frac n{2n-p}>1$ due to $p>n$.
So after finitely many steps this yields $W^{1,q}$-regularity for some 
$q=\frac{p_N}2>n$.
The same is the case if $p>2n$ at the beginning.
Next, if $u^*A$ is of class $W^{k,q}$ for some $k\in\N$, then the inhomogeneous  
terms also are of class $W^{k,q}$ and the lemma asserts the 
$W^{k+1,q}$-regularity of $\a$ and hence $u^*A$.
Iterating this argument proves the smoothness of $u^*A=\tA+\a$.
\QED

\subsection*{Weakly flat connections over a Riemann surface}

Now we consider more closely the special case when $\Si$ is a Riemann surface.
Theorem~\ref{flat regularity} shows that the injection
$\cA_{\rm flat}(\Si)/\cG(\Si) \hookrightarrow \cA^{0,p}_{\rm flat}(\Si)/\cG^{1,p}(\Si)$
in fact is a bijection.
These moduli spaces are identified and denoted by $M_\S$.
Furthermore, the holonomy induces an injection from $M_\Si$ to the
space of conjugacy classes of homomorphisms from $\pi_1(\Si)$ to $\rG$
(see e.g.\ \cite[Proposition 2.2.3]{DK}),
$$
M_\Si \,:=\;  \cA^{0,p}_{\rm flat}(\Si)/\cG^{1,p}(\Si)
\;\cong\;  \cA_{\rm flat}(\Si)/\cG(\Si)
\;\hookrightarrow\;  {\rm Hom}(\pi_1(\Si),\rG)/\sim .
$$
If $\rG$ is connected and simply connected, then every $\rG$-bundle over a Riemann surface
is automatically trivial and the holonomy in fact induces a bijection.
If there exist nontrivial $\rG$-bundles with flat connections, then
${\rm Hom}(\pi_1(\Si),\rG)/\sim$ is identified with the union of the moduli spaces for 
all such bundles.
From this one sees that $M_\Si$ is a finite dimensional singular manifold.

For $\rG={\rm SU}(2)$ for example, $M_\Si\cong{\rm Hom}(\pi_1(\Si),\SU(2))/\sim\;$
has singularities at the product connection and at the further reducible connections 
\footnote{
A connection $A\in\cA_{\rm flat}(\S)$ is called reducible if its isotropy subgroup
of $\cG(\S)$ (the group of gauge transformations that leave $A$ fixed) is not discrete.
}
-- corresponding to the connections for which the holonomy group is not ${\rm SU}(2)$ but only 
$\{\one\}$ or is conjugate to the maximal torus $S^1\subset {\rm SU}(2)$.
\footnote{
The holonomy group of a connection is given by the holonomies of all loops in $\S$.
Now the isotropy subgroup of $\cG(\S)$ of the connection is
isomorphic to the centralizer of the holonomy group, see \cite[Lemma 4.2.8]{DK}.
}
Away from these singularities, the dimension of $M_\Si$ is $6g-6$, where
$g$ is the genus of $\Si$. (The arguments in \cite[\pg 4]{DS1} show that 
$\rT_{[A]} M_\S \cong \ker\rd_A / \im \rd_A = h^1_A$ has dimension $3\cdot(2g-2)$
at irreducible connections $A$.)

For the same reasons, the space of weakly flat connections
$\cA^{0,p}_{\rm flat}(\Si)$ is in general not a Banach submanifold of 
$\cA^{0,p}(\Si)$ but a principal bundle over a singular base manifold. 
To be more precise fix a point $z\in\S$ and consider the space of based gauge 
transformations, defined as
$$
\cG^{1,p}_z(\Si) := \left\{ u\in\cG^{1,p}(\Si) \st u(z)=\one \right\} .
$$
This Lie group acts freely on $\cA^{0,p}_{\rm flat}(\Si)$.
The quotient space $\cA^{0,p}_{\rm flat}(\Si)/\cG^{1,p}_z(\Si)$ can be identified 
with ${\rm Hom}(\pi_1(\Si),\rG)$ 
(or a subset thereof if there exist nontrivial $\rG$-bundles over $\S$) 
via the holonomy based at $z$.
This based holonomy map $\r_z:\cA^{0,p}_{\rm flat}(\Si)\to{\rm Hom}(\pi_1(\Si),\rG)$ is
defined by first choosing a based gauge transformation that makes the 
connection smooth and then computing the holonomy around loops based at $z$.
Now $\r_z$ gives $\cA^{0,p}_{\rm flat}(\Si)$ the structure of a principal bundle 
with fibre $\cG^{1,p}_z(\Si)$ over the finite dimensional singular manifold 
${\rm Hom}(\pi_1(\Si),\rG)$ (or a subset thereof)
$$
\cG^{1,p}_z(\Si) 
\hookrightarrow  \cA^{0,p}_{\rm flat}(\Si)
\overset{\r_z}{\longrightarrow} {\rm Hom}(\pi_1(\Si),\rG) .
$$
Note that this discussion does not require the Riemann surface $\Si$ to be connected.
Only when fixing a base point for the holonomy map and the based
gauge transformations one has to adapt the definition.
Whenever $\Si=\bigcup_{i=1}^n \Si_i$ has several connected components $\Si_i$, then
'fixing a point $z\in \Si$' implicitly means that one fixes a point $z_i\in \Si_i$
in each connected component. The group of based gauge transformations then
becomes
$$
\cG^{1,p}_z({\textstyle\bigcup_{i=1}^n \Si_i}) 
:= \left\{ u\in\cG^{1,p}(\Si) \st u(z_i)=\one \quad\forall i=1,\ldots,n \right\} .
$$

\pagebreak

\section{Lagrangians in the space of connections}
\label{spaces}

Consider the trivial $\rG$-bundle over a (possibly disconnected) Riemann surface $\Si$ 
of (total) genus $g$, where $\rG$ is a compact Lie group with Lie algebra $\cg$.
There is a gauge invariant symplectic form $\o$ on the space of connections
$\cA^{0,p}(\Si)$ for $p>2$ defined as follows.
For tangent vectors $\a,\b\in L^p(\Si,\rT^*\Si\otimes\cg)$ to the affine space $\cA^{0,p}(\Si)$
\begin{equation} \label{symp form}
\o(\a,\b) = \int_\Si \la \a\wedge\b \ra .
\end{equation}
The action of the infinite dimensional gauge group $\cG^{1,p}(\Si)$ on the 
symplectic Banach space $(\cA^{0,p}(\Si),\o)$ is Hamiltonian with moment map 
$A\mapsto *F_A$ (more precisely, the equivalent weak expression in 
$(W^{1,p^*}(\Si,\cg))^*$).
So the moduli space of flat connections $M_\S=\cA^{0,p}_{\rm flat}(\Si)/\cG^{1,p}(\Si)$ 
can be viewed as the symplectic quotient $\cA^{0,p}(\Si)/\hspace{-1mm}/\cG^{1,p}(\Si)$
as was first observed by Atiyah and Bott \cite{AB}.
However, $0$ is not a regular value of the moment map, so $M_\Si$ is a singular 
symplectic manifold. Due to these singularities at the reducible connections the 
infinite dimensional setting suggests itself.

Note that for any metric on $\S$ the Hodge $*$ operator is an $\o$-compatible complex structure 
since $\o(\cdot,*\cdot)$ is the $L^2$-metric:
For all $\a,\b\in L^p(\Si,\rT^*\Si\otimes\cg)$
\begin{equation}\label{om,L2}
\o(\a,*\b) \;=\; \int_\S \la \a \wedge *\b \ra \;=\; \la \a \,,\, \b \ra_{L^2}.
\end{equation}
Next, we call a Banach submanifold $\cL\subset\cA^{0,p}(\Si)$ Lagrangian if it is isotropic, 
i.e.\ $\o|_\cL \equiv 0$, and if $\rT_A\cL$ is maximal for all $A\in\cL$ in the following sense:
If $\o(\rT_A\cL,\a)=\{0\}$ for some $\a\in\cA^{0,p}(\Si)$, then $\a\in\rT_A\cL$.
In general, this condition does not imply that $\cL$ is also totally real with respect to any
$\o$-compatible complex structure. 
However, we will only consider Lagrangian submanifolds $\cL\subset\cA^{0,p}(\S)$ that are gauge 
invariant and contained in the space of weakly flat connections. These are automatically totally real
with respect to the Hodge $*$ operator, as lemma \ref{totally real} will show. 
It is based on the following twisted Hodge decomposition.

\begin{lem}  \label{HD}
Fix a metric on $\S$ and let $A\in\cA^{0,p}_{\rm flat}(\S)$. Then 
$$
L^p(\S,\rT^*\S\otimes\cg) 
\;=\; \rd_A W^{1,p}(\S,\cg) \,\oplus\, *\rd_A W^{1,p}(\S,\cg) \,\oplus\, h^1_A ,
$$
with the finite dimensional space 
$h^1_A=\ker \rd_A \cap \ker \rd_A^* \subset W^{1,q}(\S,\rT^*\S\otimes\cg)$
for $\frac 1q = \frac 12 + \frac 1p$.
\end{lem}
\Pr
Recall that $p>2$, hence $\rd_A W^{1,2}(\S,\cg)\subset L^2(\S,\rT^*\S\otimes\cg)$
due to the Sobolev embedding $W^{1,2}(\S)\hookrightarrow L^r(\S)$ for any $r<\infty$.
The weak flatness of $A$ then implies that $\rd_A W^{1,2}(\S,\cg)$ and $*\rd_A W^{1,2}(\S,\cg)$
are $L^2$-orthogonal. The orthogonal complement of their direct sum then exactly is $h^1_A$.
(To see that every $L^2$-regular $1$-form that is orthogonal to $\im\rd_A$ and to $*\im\rd_A$ is 
automatically $W^{1,q}$-regular, one can use the regularity theory for the Laplace operator or the Hodge 
decomposition, or see e.g.\ \cite[Lemma A.2]{W ell}.)
Next, note that \hbox{$h^1_A\subset L^p(\S,\rT^*\S\otimes\cg)$} due to the Sobolev embedding
$W^{1,q}(\S)\hookrightarrow L^p(\S)$. Now the same regularity arguments as above show that
the orthogonal decomposition 
\begin{equation} \label{HD 2}
L^2(\S,\rT^*\S\otimes\cg) 
\;=\; \rd_A W^{1,2}(\S,\cg) \,\oplus\, *\rd_A W^{1,2}(\S,\cg) \,\oplus\, h^1_A 
\end{equation}
restricts to the claimed decomposition of $L^p(\S,\rT^*\S\otimes\cg)$.
Finally, to see that $h^1_A$ is finite dimensional note that it is isomorphic to the cokernel
of the operator $\rd_A\oplus*\rd_A: W^{1,p}(\S,\cg)\times W^{1,p}(\S,\cg) \to L^p(\S,\rT^*\S\otimes\cg)$.
Now this operator is a compact perturbation of the Fredholm operator $\rd\oplus*\rd$, hence its
cokernel is finite dimensional.
\QED

\begin{lem}  \label{totally real}
Let $\cL\subset\cA^{0,p}(\S)$ be a Lagrangian submanifold.
Suppose that $\cL\subset\cA^{0,p}_{\rm flat}(\S)$ and that $\cL$ is invariant 
under the action of $\cG^{1,p}(\S)$.
Then $\cL$ is totally real with respect to the Hodge $*$ operator for any
metric on $\S$.
\end{lem}
\Pr
Pick any $A\in\cL$ and denote $L:=\rT_A\cL$. Then we have to show that $\cA^{0,p}(\S)=L\oplus *L$.
Firstly, the only element $\a\in L\cap *L$ in the intersection is $\a=0$ since $*\a \in **L=L$ 
and thus $\|\a\|_{L^2}^2 = \o(\a,*\a) = 0$.

Secondly, to see that the direct sum $L\oplus *L$ exhausts all of $\cA^{0,p}(\S)$, assume the
contrary. Then there exists a nonzero linear functional $\p$ on $\cA^{0,p}(\S)$ that vanishes
on $L\oplus *L$. Due to the gauge invariance of $\cL$ one has $\rd_A W^{1,p}(\S,\cg)\subset L$,
so $\p$ vanishes in particular on $ \rd_A W^{1,p}(\S,\cg) \,\oplus\, *\rd_A W^{1,p}(\S,\cg)$.
Now recall the Hodge decomposition in lemma \ref{HD} and (\ref{HD 2}). It implies that $\p$ has to be 
nonzero on $h^1_A$ and hence can be extended to a nonzero linear functional on $\cA^{0,2}(\S)$ that 
vanishes on $ \rd_A W^{1,2}(\S,\cg) \,\oplus\, *\rd_A W^{1,2}(\S,\cg)$.
Thus the extended functional can be written as $\p=\la \a, \cdot \ra_{L^2}$ for some 
$\a\in L^2(\S,\rT^*\S\otimes\cg)$.
But now the orthogonal decomposition (\ref{HD 2}) implies that $\a\in h^1_A\subset\cA^{0,p}(\S)$. 
Now for all $\b\in L=\rT_A\cL$ one has
$$
\o(\b,\a) \;=\; \la \a , *\b \ra_{L^2} \;=\; \p(*\b) \;=\; 0 .
$$
The Lagrangian property of $\cL$ then implies that $\a\in L$ and hence
$$
\|\a\|_{L^2}^2 \;=\; \la \a,\a \ra_{L^2} \;=\; \p(\a) \;=\; 0 .
$$
This proves $\a=0$ in contradiction to the assumption $\p\neq 0$. Hence $\cL$ is indeed
totally real with respect to the complex structure $*$, i.e.\ for all $A\in\cL$
\begin{equation}\label{L}
L^p(\S,\rT^*\S\otimes\cg) \;=\; \rT_A\cL \,\oplus\, *\,\rT_A\cL .
\end{equation}
\QED

The assumption $\cL\subset\cA^{0,p}_{\rm flat}(\Si)$ directly implies that $\cL$ is
gauge invariant if $\rG$ is connected and simply connected. On the other hand, 
the gauge invariance of $\cL$ implies $\cL\subset\cA^{0,p}_{\rm flat}(\Si)$ if the 
Lie bracket on $\rG$ is nondegenerate (i.e.\ the center of $\rG$ is discrete).
So for example in the case $\rG={\rm SU}(2)$ both conditions are equivalent.
We will always assume both conditions.
Then moreover, $\cL$ descends to a (singular) submanifold of the (singular) moduli space
of flat connections,
$$
L \,:=\; \cL/\cG^{1,p}(\Si) \subset \cA_{\rm flat}^{0,p}(\Si)/\cG^{1,p}(\Si) \;=:\, M_\Si .
$$
This submanifold is obviously isotropic, i.e.\ the symplectic structure induced by
(\ref{symp form}) on $M_\Si$ vanishes on $L$.
Moreover, its tangent spaces have half of the dimension of those of $M_\Si$, so 
$L\subset M_\Si$ is a Lagrangian submanifold.
Indeed, in the Hodge decomposition, lemma \ref{HD}, $*\rd_A W^{1,p}(\S,\cg)$ is the complement of
$\ker \rd_A = \rT_A\cA_{\rm flat}^{0,p}(\Si)$, $\rd_A W^{1,p}(\S,\cg)$ is the tangent space to the 
orbit of $\cG^{1,p}(\Si)$ through $A$, and so $h^1_A\cong\rT_{[A]}M_\Si$.
Now compare this with the decomposition (\ref{L}). Here $\rT_A\cL=\rd_A W^{1,p}(\S,\cg)\oplus V$, 
where the complement $V\subset\cA^{0,p}(\S)$ is finite dimensional and $V\oplus *V$ can replace 
$h^1_A$ in the Hodge decomposition.
Thus $\rT_{[A]}\cL \cong \rT_A\cL/\rd_A W^{1,p}(\S,\cg) \cong V$ must have half the dimension 
of $h^1_A$.

Moreover, our assumptions on the Lagrangian submanifold ensure that the holonomy map 
$\r_z:\cL\to {\rm Hom}(\pi_1(\Si),\rG)$ based at $z\in\Si$ is welldefined and invariant 
under the action of the based gauge group $\cG^{1,p}_z(\Si)$.
(The holonomy map and based gauge group are introduced in section~\ref{weakly flat}.)
Note that ${\rm Hom}(\pi_1(\Si),\rG)$ naturally embeds into
${\rm Hom}(\pi_1(\Si\setminus\{z\}),\rG)$, which is a smooth manifold 
diffeomorphic to $G^{2g}$. This gives ${\rm Hom}(\pi_1(\Si),\rG)$ a 
differentiable structure (that is in fact independent of $z\in\Si$),
however, it is a manifold with singularities.
In the following lemma we list some crucial properties of the Lagrangian submanifolds. 
Here we use the notation
$$
W^{1,p}_z(\Si,\cg) :=  \bigl\{ \x\in W^{1,p}(\Si,\cg) \st \x(z)=0 \bigr\} 
$$
for the Lie algebra $\rT_\smone\cG^{1,p}_z(\S)$ of the based gauge group.
(If $\Si$ is not connected then as before one fixes a base point in each 
connected component and modifies the definition of $W^{1,p}_z(\Si,\cg)$ 
accordingly.)
Moreover, we will denote the differential of a map $\p$ at a point $x$ by $\rT_x\p$ in order 
to distinguish it from the exterior differential on differential forms, $\rd_A$, associated 
with a connection $A$.

\begin{lem}  \label{Lagrangian lemma}
Let $\cL\subset\cA^{0,p}(\Si)$ be a Lagrangian submanifold and fix $z\in\Si$. 
Suppose that $\cL\subset\cA^{0,p}_{\rm flat}(\Si)$ and that $\cL$ is invariant 
under the action of $\cG^{1,p}(\Si)$. Then the following holds:
\begin{enumerate}

\item  $L:=\cL/\cG^{1,p}_z(\Si)$ is a smooth manifold of dimension 
$m=g\cdot\dim \rG$ and the holonomy induces a diffeomorphism $\r_z:L\to M$ to a 
submanifold $M\subset{\rm Hom}(\pi_1(\Si),\rG)$.

\item  $\cL$ has the structure of a principal $\cG^{1,p}_z(\Si)$-bundle over 
$M$,
$$
\cG^{1,p}_z(\Si) 
\hookrightarrow  \cL \overset{\r_z}{\longrightarrow} M .
$$

\item  Fix $A\in\cL$. Then there exists a local section $\p:V \to \cL$ over a 
neighbourhood $V\subset\R^m$ of \,$0$ such that $\p(0)=A$ and $\r_z\comp\p$ is a 
diffeomorphism to a neighbourhood of $\r_z(A)$.
This gives rise to Banach submanifold coordinates for $\cL\subset\cA^{0,p}(\Si)$,
namely a smooth embedding
$$
\Th: \cW \to \cA^{0,p}(\Si)
$$
defined on a neighbourhood
$\cW\subset W^{1,p}_z(\Si,\cg) \times \R^m \times W^{1,p}_z(\Si,\cg) \times \R^m$
of zero by
$$
\Th(\x_0,v_0,\x_1,v_1) := \exp(\x_0)^*\p(v_0) + *\rd_A\x_1 + *\rT_0\p(v_1) .
$$
Moreover, if $A$ is smooth, then the local section can be chosen such that 
the image is smooth, $\p:V \to \cL\cap\cA(\Si)$. Now the same map $\Th$ is a diffeomorphism
between neighbourhoods of zero in $W^{1,p}_z(\Si,\cg^2) \times \R^{2m}$
and neighbourhoods of $A$ in $\cA^{0,p}(\Si)$ for all $p>2$.
\end{enumerate}
\end{lem}

We postpone the proof and first note that this lemma shows that the Lagrangian submanifolds 
considered here all satisfy the crucial assumption for theorem~\ref{holo regularity} and 
\ref{holo regularity 2}.

\begin{cor}  \label{Hp}
Let $\cL\subset\cA^{0,p}(\Si)$ be as in lemma \ref{Lagrangian lemma}, then it satisfies
$(H_p)$, i.e.\ $\cL$ is modelled on a closed subspace of an $L^p$-space.
\end{cor}

\noindent
{\bf Proof of corollary \ref{Hp}: } \\
The bundle structure of $\cL$ in lemma~\ref{Lagrangian lemma}~(ii) shows that $\cL$ is 
modelled on $W^{1,p}_z(\S,\cg)\times\R^m$.
This is since the Banach manifold $\cG^{1,p}_z(\S)$ is modelled on $W^{1,p}_z(\S,\cg)$, 
which is a closed subspace of $W^{1,p}(\S,\cg)$ and thus is isomorphic to a closed subspace 
of $L^p(\S,\cg^3)$.
Now recall example~\ref{Hp ex} to see that $\cL$ indeed satisfies $(H_p)$.
\QED

The Banach submanifold charts $\Th$ in lemma~\ref{Lagrangian lemma}~(iii) are essentially the 
same as the charts $\Th$ in the proof of theorem~\ref{holo regularity}.
In this special case, we have more detailed information on the structure of $\Th$, which is
the main point in the proof of the following approximation result for $W^{1,p}$-connections 
with Lagrangian boundary values.

\begin{cor}  \label{bc approx}
Let $\cL\subset\cA^{0,p}(\Si)$ be as in lemma \ref{Lagrangian lemma} and let
$$
\Om\subset\H:=\{(s,t)\in\R^2 \st t\geq 0\}
$$ 
be a compact submanifold.
Suppose that $A\in \cA^{1,p}(\Om\times\Si)$ satisfies the boundary 
condition
\begin{equation}\label{bc}
A|_{(s,0)\times\Si} \in \cL 
\qquad\forall (s,0)\in\pd\Om .
\end{equation}
Then there exists a sequence of smooth connections $A^\n\in\cA(\Om\times\Si)$
that satisfy (\ref{bc}) and converge to $A$ in the $W^{1,p}$-norm.
\end{cor}

\noindent
{\bf Proof of corollary \ref{bc approx}: } \\
We decompose $A=\P\ds+\Psi\dt+B$ into two functions 
$\P,\Psi\in W^{1,p}(\Om\times\Si,\cg)$ and a family of $1$-forms 
$B\in W^{1,p}(\Om\times\Si,\rT^*\Si\otimes\cg)$ on $\Si$ such that 
$B(s,0)\in\cL$ for all $(s,0)\in\pd\Om$.
Then it suffices to find an approximating sequence for $B$ with Lagrangian 
boundary conditions on a neighbourhood of $\Om\cap\pd\H$. This can be patched 
together with any smooth $W^{1,p}$-approximation of $B$ on the rest of $\Om$
and can be combined with standard approximations of the functions $\P$ and 
$\Psi$ to obtain the required approximation of $A$.

So fix any $(s_0,0)\in\Om\cap\pd\H$ and use theorem~\ref{flat regularity} to 
find $u_0\in\cG^{1,p}(\Si)$ such that $A_0:=u_0^*B(s_0,0)$ is smooth.
Lemma~\ref{Lagrangian lemma}~(iii) gives a \hbox{diffeomorphism}
$\Th:\cW \to\cV$ between neighbourhoods 
$\cW\subset W^{1,p}_z(\Si,\cg^2)\times\R^{2m}$ of zero and
$\cV\subset\cA^{0,p}(\Si)$ of $A_0$. This was constructed such that 
$\cC^\infty(\Si,\cg^2)\times\R^{2m}$ is mapped to $\cA(\Si)$ and such that
$\Th:\cW\cap W^{2,p}_z(\Si,\cg^2)\times\R^{2m} \to \cV\cap\cA^{1,p}(\Si)$ 
also is a diffeomorphism.
Now note that $B\in\cC^0(\Om,\cA^{0,p}(\Si))$.
Hence there exists a neighbourhood $U\subset\Om$ of $(s_0,0)$ and one can choose a smooth gauge 
transformation $u\in\cG(\Si)$ that is $W^{1,p}$-close to $u_0$ such that
$u^*B(s,t)\in\cV$ for all $(s,t)\in U$.
Now we define $\x=(\x_0,\x_1):U\to W^{1,p}_z(\Si,\cg^2)$ and $v=(v_0,v_1):U\to\R^{2m}$ by
$\Th(\x(s,t),v(s,t))=u^*B(s,t)$. 
Recall that $B$ is of class $W^{1,p}$ on $U\times\S$, hence it lies in both
$W^{1,p}(U,\cA^{0,p}(\S)$ and $L^p(U,\cA^{1,p}(\S)$.
Thus $\x \in W^{1,p}(U,W^{1,p}_z(\Si,\cg^2)) \cap L^p(U,W^{2,p}_z(\Si,\cg^2))$
and $v\in W^{1,p}(U,\R^{2m})$, and these satisfy the boundary conditions $\x_1|_{t=0}=0$ and 
$v_1|_{t=0}=0$ due to the Lagrangian boundary condition for $B$. 
Now there exist $\x^\n\in\cC^\infty(U\times\Si,\cg^2)$ and $v^\n\in\cC^\infty(U,\R^{2m})$ 
such that $\x^\n\to\x$ and $v^\n\to v$ in all these spaces,
$\x^\n(\cdot,z)\equiv 0$, $\x^\n_1|_{t=0}=0$, and $v^\n_1|_{t=0}=0$.
(These are constructed with the help of mollifiers as in lemma~\ref{approx}. 
One first reflects $\x$ at the boundary and mollifies it with
respect to $U$ to obtain approximations in $\cC^\infty(U,W^{2,p}_z(\S,\cg^2))$ with zero
boundary values. Next, one mollifies on $\S$, and finally one corrects the value at $z$.)
It follows that $B^\n(s,t):=(u^{-1})^*\Th(\x(s,t),v(s,t))$ is a sequence of
smooth maps from $U$ to $\cA(\Si)$ which satisfies the Lagrangian boundary 
condition and converges to $B$ in the $W^{1,p}$-norm.

Now $\Om\cap\pd\H$ is compact, so it is covered by finitely many such 
neighbourhoods $U_i$ on which there exist smooth $W^{1,p}$-approximations of $B$
with Lagrangian boundary values.
These can be patched together in a finite procedure since the above construction 
allows to interpolate in the coordinates between $\x^\n$, $v^\n$ and other smooth
approximations $\x'$, $v'$ (arising from approximations of $B$ on another 
neighbourhood $U'$ in different coordinates) of $\x$ and $v$ respectively.
This gives the required approximation of $B$ in a neighbourhood of $\Om\cap\pd\H$ 
and thus finishes the proof.
\QED \\

\noindent
{\bf Proof of lemma \ref{Lagrangian lemma}: } \\
Fix $A\in\cL$ and consider the following two decompositions:
\begin{align}
L^p(\Si,\rT^*\Si\otimes\cg) 
&= \rT_A\cL \oplus *\rT_A\cL  \label{decomp}\\
&= \rd_A W^{1,p}_z(\Si,\cg) \oplus *\rd_A W^{1,p}_z(\Si,\cg) \oplus \tilde h_A .
\nonumber
\end{align}
The first direct sum is due to lemma \ref{totally real}.
In the second decomposition, $\tilde h_A$ is a complement of the image of the 
following Fredholm operator:
$$
D_A : \begin{array}{ccc}
W^{1,p}_z(\Si,\cg) \times W^{1,p}_z(\Si,\cg) &\longrightarrow
                                           & L^p(\Si,\rT^*\Si\otimes\cg) \\
(\x,\z) &\longmapsto& \rd_A\x + * \rd_A\z .
\end{array}
$$
To see that $D_A$ is Fredholm note that for every $A\in\cA^{0,p}_{\rm flat}(\Si)$ the operator 
$D_A$ is injective and is a compact perturbation of $D_0$.
Hence the dimension of ${\rm coker}\,D_A$ (and thus of $\tilde h_A$) is the same
as that of ${\rm coker}\,D_0$. In the case $A=0$ one can choose the space of $\cg$-valued
harmonic $1$-forms $h^1=\ker\rd\cap\ker\rd^*$ as complement $\tilde h_0$. 
So $\tilde h_A$ must always have the dimension 
${\rm dim}\,\tilde h_A ={\rm dim}\,h^1 = 2g \cdot {\rm dim}\,\rG =2m$.
(Note that in general one can choose $\tilde h_A$ to contain $h^1_A$, 
but this might not exhaust the whole complement.)

Due to the $\cG^{1,p}_z(\Si)$-invariance of $\cL$ the splittings (\ref{decomp})
now imply that there exists an $m$-dimensional subspace $L_A\subset \tilde h_A$
such that
$$
\rT_A\cL = \rd_A W^{1,p}_z(\Si,\cg) \oplus L_A .
$$
So $\rT_A\cL$ is isomorphic to the Banach space $W^{1,p}_z(\Si,\cg)\times\R^m$
via $\rd_A\oplus F$ for some isomorphism $F:\R^m\to L_A$.
Here we have used the fact that $\rd_A$ is injective when restricted to $W^{1,p}_z(\Si,\cg)$.
Now choose a coordinate chart $\P:\rT_A\cL \to \cL$ defined near $\P(0)=A$, 
then the following map is defined for a sufficiently small neighbourhood 
$V\subset\R^m$ of $0$,
$$
\Psi:   \begin{array}{ccc}
\cG^{1,p}_z(\Si) \times V &\longrightarrow & \cL \\
(u,v) &\longmapsto& u^*\bigl( \P\comp(\rT_0\P)^{-1}\comp F (v) \bigr).
\end{array}
$$
We will show that this is an embedding and a submersion (and thus a diffeomorphism
to its image).
Firstly, $\rT_{(\smone,0)}\Psi : (\x,w)\mapsto \rd_A\x + F w$ is an
isomorphism. Next, note that $\Psi(u,v)=u^*\Psi(\one,v)$ and use this to 
calculate for all $u\in\cG^{1,p}_z(\Si)$, $\x\in W^{1,p}_z(\Si,\cg)$, and 
$v,w\in\R^m$
$$ 
\rT_{(u,v)}\Psi : (\x u,w) 
\mapsto u^{-1} \bigl( \rd_{\Psi(\smone,v)}\x 
                    + \rT_{(\smone,v)}\Psi (0,w) \bigr) u.
$$ 
One sees that $u(\rT_{(u,v)}\Psi ) u^{-1}$ is a small perturbation of 
$\rT_{(\smone,0)}\Psi$, hence one can choose $V$ 
sufficiently small (independently of $u$) such that $\rT_{(u,v)}\Psi$ also 
is an isomorphism for all $v\in V$.
So it remains to check that $\Psi$ in fact is globally injective.

Suppose that $u,u'\in\cG^{1,p}_z(\Si)$ and $v,v'\in V$ such that
$\Psi(u,v)=\Psi(u',v')$. Rewrite this as $\Psi(\one,v)=\Psi(\tilde u,v')$ with 
$\tilde u := u'u^{-1}\in\cG^{1,p}_z(\Si)$. 
Now by the choice of a sufficiently small $V$ the norm $\|\Psi(\one,v)-\Psi(\one,v')\|_p$ 
can be made arbitrarily small. Then the identity $\Psi(\one,v)=\tilde u^*\Psi(\one,v')$
automatically implies that $\tilde u$ is $\cC^0$-close to $\one$.
(Otherwise one would find a sequence of $L^p$-connections $A^\n\to A$ and 
$u^\n\in\cG^{1,p}_z(\Si)$ such that $\| u^{\n\;*}A^\n-A^\n \|_p \to 0$ 
but $d_{\cC^0}(u^\n,\one)\geq\D>0$. However, from
$(u^\n)^{-1}\rd u^\n = u^{\n\;*}A^\n - (u^\n)^{-1}A^\n u^\n$ one obtains an $L^p$-bound on 
$\rd u^\n$ and thus finds a weakly $W^{1,p}$-convergent subsequence of the $u^\n$.
Its limit $u\in\cG^{1,p}_z(\Si)$ would have to satisfy $u^*A=A$, hence
$u\equiv\one$ in contradiction to $d_{\cC^0}(u,\one)\geq\D>0$.)
So one can write $\tilde u=\exp(\x)$ where $\x\in W^{1,p}_z(\S,\cg)$ is small in the
$L^\infty$-norm. Next, the identity
$$
\tilde u ^{-1} \rd \tilde u 
= \Psi(\one,v) - \tilde u ^{-1}\Psi(\one,v') \tilde u
$$
shows that $\|\x\|_{W^{1,p}}$ will be small if $V$ is small (and thus $\tilde u$ is 
$\cC^0$-close to $\one$).
Hence if $V$ is sufficiently small, then $(\tilde u,v')$ and $(\one,v)$ automatically
lie in a neighbourhood of $(\one,0)$ on which $\Psi$ is injective, and hence
$u=u'$ and $v=v'$.

We have thus shown that $\Psi:\cG^{1,p}_z(\Si) \times V \to \cL$ is a 
diffeomorphism to its image.
This provides manifold charts $\psi:V\to\cL/\cG^{1,p}_z(\Si)$,
$v\mapsto [\Psi(\one,v)]$ for $L:=\cL/\cG^{1,p}_z(\Si)$.
Now fix $2g$ generators of the fundamental group $\pi_1(\Si)$ based at $z$, then the 
corresponding holonomy map $\r_z:L\to \rG\times\cdots\times \rG$ is an embedding, 
so its image $M\subset{\rm Hom}(\pi_1(\Si),\rG)$ is a smooth submanifold.
This proves (i).
For (ii) the diffeomorphism $\Psi$ gives a bundle chart over $\cU:=\r_z(\psi(V))\subset M$,
namely 
$$
\Psi\comp\bigl({\rm id}\times (\r_z\comp\psi)^{-1} \bigr) \,:\; 
\cG^{1,p}_z(\Si) \times \cU \longrightarrow \cL .
$$
Furthermore, the local section for (iii) is given by $\p(v):=\Psi(\one,v)$.
However, this is a map $\p:V\to\cL$; it does not necessarily take values in
the smooth connections.
Now if $A\in\cL\cap\cA(\Si)$ is smooth, then for a sufficiently small neighbourhood $V$
this section can be modified by gauge transformations such that $\p:V\to\cL\cap\cA(\Si)$.
To see this, note that the gauge transformations in the local slice theorem are given
by an implicit function theorem: One solves $D(v,\x)=0$ for $\x=\x(v)\in W^{1,p}(\S,\cg)$
with the following operator:
$$
D :
\begin{array}{cccl}
V \times W^{1,p}(\S,\cg) &\longrightarrow & \im\rd'_A &\subset 
\;\bigl( W^{1,p^*}(\S,\cg)\bigr)^*\\
(v,\x) &\longmapsto& \rd'_A \bigl(\exp(\x)^*\p(v) - A\bigr) &.
\end{array}
$$
Here $\rd'_A$ denotes the dual operator of $\rd_A$ on $W^{1,p^*}(\S,\cg)$.
One has \hbox{$D(0,0)=0$} and checks that $\pd_2 D(0,0):\x\to \rd'_A\rd_A\x$ is a surjective map 
to $\im\rd'_A$, see e.g.\ \cite[Lemma 8.5]{W}.
The implicit function theorem \cite[XIV,Theorem 2.1]{Lang} then gives the required gauge
transformations $\exp(\x(v))\in\cG^{1,p}(\S)$ that bring $\p(v)$ into local Coulomb gauge
and thus make it smooth. (By construction $\p(v)$ is weakly flat, then see the proof of 
theorem~\ref{flat regularity}.)
This modification by gauge transformations however does not affect the topological 
direct sum decomposition $\rT_A\cL = \rd_A W^{1,p}_z(\Si,\cg) \oplus \im\rT_0\p$.

To see that the given map $\Th$ is a diffeomorphism between 
neighbourhoods of $0$ and $A$ just note that the inverse of $\rT_0\Th$ is 
given by the splitting
\begin{align*}
L^p(\Si,\rT^*\Si\otimes\cg) 
&= \rT_A\cL \oplus *\rT_A\cL  \\
&= \rd_A W^{1,p}_z(\Si,\cg) \oplus \im\rT_0\p
   \oplus *\rd_A W^{1,p}_z(\Si,\cg) \oplus *\im\rT_0\p
\end{align*}
composed with the inverses of $\rd_A|_{W^{1,p}_z(\Si,\cg)}$ and $\rT_0\p$.
\QED

Now observe that the choice of $p>2$ for the Lagrangian submanifolds in the 
above lemma is accidental. All connections $A\in\cL$ are gauge equivalent to a 
smooth connection, and the $L^q$-completion of $\cL\cap\cA(\Si)$ is a Lagrangian 
submanifold in $\cA^{0,q}(\Si)$ for all $q>2$.
In fact, this simply is the restricted ($q>p$) or completed ($q<p$)
$\cG^{1,q}_z(\Si)$-bundle over $M$.

\subsection*{The main example}

Suppose that $\rG$ is connected and simply connected and that $\Si=\pd Y$ is the 
boundary of a handlebody $Y$. 
(Again, the handlebody and thus its boundary might consist of several connected 
components.)
This together with the fact that $\pi_2(\rG)=0$ for any Lie group $\rG$ (see e.g.\
\cite[Proposition 7.5]{Broecker}) ensures that the gauge group $\cG^{1,p}(\Si)$ is 
connected and that every gauge transformation on $\Si$ can be extended to $Y$.

Let $p>2$ and let $\cL_Y$ be the $L^p(\Si)$-closure of the set of smooth flat 
connections on $\Si$ that can be extended to a flat connection on $Y$,
$$
\cL_Y \,:=\; {\rm cl}\, \bigl\{ A\in\cA_{\rm flat}(\Si) \st \exists 
                              \tA\in\cA_{\rm flat}(Y) : \tA|_\Si=A \bigr\} 
\;\subset\;\cA^{0,p}(\Si).
$$
This is an example of a totally real submanifold of $(\cA^{0,p}(\S,\cg),*)$ that satisfies
the assumption of theorem~\ref{holo regularity} and \ref{holo regularity 2}. 
This is due to the lemmata~\ref{totally real} and \ref{Lagrangian lemma} and the 
following properties of $\cL_Y$.


\begin{lem} \hspace{1mm} \label{LY} \\
\vspace{-5mm}
\begin{enumerate}
\item
$\cL_Y = 
 \bigl\{ u^*(A|_\Si) \st A\in\cA_{\rm flat}(Y), u\in\cG^{1,p}(\Si) \bigr\} $
\item  $\cL_Y\subset\cA^{0,p}(\Si)$ is a Lagrangian submanifold.
\item  $\cL_Y\subset\cA^{0,p}_{\rm flat}(\Si)$ and $\cL_Y$ is invariant under 
       the action of $\cG^{1,p}(\Si)$.
\item  Fix any $z\in\Si$. Then
$$
\cL_Y = \bigl\{ A\in\cA^{0,p}_{\rm flat}(\Si) \st 
    \r_z(A) \in {\rm Hom}(\pi_1(Y),\rG) \subset {\rm Hom}(\pi_1(\Si),\rG)  \bigr\} .
$$
So $\cL_Y$ obtains the structure of a $\cG^{1,p}_z(\Si)$-bundle 
over the $g$-fold product 
$M=\rG\times \cdots \times \rG \cong {\rm Hom}(\pi_1(Y),\rG)$,
$$
\cG^{1,p}_z(\Si) 
\hookrightarrow  \cL_Y
\overset{\r_z}{\longrightarrow} {\rm Hom}(\pi_1(Y),\rG) .
$$
\end{enumerate}
\end{lem}
\Pr
Firstly, $\cL_Y\subset\cA^{0,p}_{\rm flat}(\Si)$ follows from the fact that 
weak flatness is an $L^p$-closed condition for $p> 2$.
The holonomy $\r_z:\cA^{0,p}_{\rm flat}(\Si)\to \rG\times\cdots\times \rG$ 
is continuous with respect to the $L^p$-topology. Thus for every $A\in\cL_Y$
the holonomy vanishes on those loops in $\Si$ that are contractible in $Y$. 
On the other hand, in view of theorem~\ref{flat regularity}, every 
$A\in\cA^{0,p}_{\rm flat}(\Si)$ whose holonomy descends to 
${\rm Hom}(\pi_1(Y),\rG)$ can be written as $A=u^*\tA$, where 
$u\in\cG^{1,p}_z(\Si)$ and the holonomy of $\tA\in\cA_{\rm flat}(\Si)$ also 
vanishes along the loops that are contractible in $Y$.
Thus $\tA$ can be extended to a flat connection on $Y$ and smooth approximation
of $u$ proves that $A\in\cL_Y$. This proves the alternative definitions of 
$\cL_Y$ in (iv) and (i). Then (iii) is a consequence of (i).

To prove the second assertion in (iv) we explicitly construct local sections of $\cL_Y$.
Let the loops $\a_1,\b_1,\ldots,\a_g,\b_g \subset \Si$ be disjoint from $z$ and
represent the standard 
generators of $\pi_1(\Si)$ such that $\a_1,\ldots,\a_g$ generate $\pi_1(Y)$ and 
such that the only nonzero intersections are $\a_i \cap \b_i$.
One can then modify the $\a_i$ such that they run through $z$ but still do not 
intersect the $\b_j$ for $j\neq i$.
\footnote{$\pi_1(\Si)$ is the quotient of the free group generated by elements
$\a_1,\b_1,\ldots,\a_g,\b_g$ by the relation 
$\a_1\b_1\a_1^{-1}\b_1^{-1}\ldots\a_g\b_g\a_g^{-1}\b_g^{-1}=\one$,
whereas $\pi_1(Y)$ is the free group generated by $\a_1,\ldots,\a_g$.
}
Now fix $A\in\cL_Y$.
In order to change its holonomy along $\a_i$ by some $g\in \rG$ close to $\one$, 
one gauge transforms $A$ in a small neighbourhood of $\b_i$ in $\S$ with a smooth gauge 
transformation that equals $\one$ and $g$ respectively near the two boundary 
components of that ring about $\b_i$.
That way one obtains a smooth local section $\p : V \to \cL_Y$ defined on a 
neighbourhood $V\subset \cg^g$ of $0$, such that $\p(0)=A$ and
\hbox{$\r_z\comp\p : V \to {\rm Hom}(\pi_1(Y),\rG)$} is a bijection onto a neighbourhood 
of $\r_z(A)$.
This leads to a bundle chart
$$
\Psi: 
\begin{array}{ccc}
\cG^{1,p}_z(\Si) \times V &\longrightarrow& \cL_Y \\
(u,v) &\longmapsto& u^*\p(v) .
\end{array}
$$
%
%
Note that for smooth $A\in\cL_Y\cap\cA(\Si)$ the local section $\p$ constructed
above in fact is a section in the smooth part $\cL_Y\cap\cA(\Si)$ of the
Lagrangian.
Using these bundle charts one also checks that $\cL_Y\subset\cA^{0,p}(\Si)$ is 
indeed a Banach submanifold. Now a submanifold chart near $\Psi(u,v)\in\cA^{0,p}(\S)$ is 
given by $(\x,w)\mapsto \Psi(\exp(\x)u,v+w) + *\rT_{(u,v)}\Psi(\x,w)$. 
As in lemma~\ref{Lagrangian lemma} one checks that this is a local diffeomorphism.

To verify the Lagrangian condition it suffices to 
consider $\o$ on $\rT_A\cA^{0,p}(\Si)$ for smooth $A\in\cL_Y$.
This is because both $\o$ and $\cL_Y$ are invariant under the gauge action.
So we pick some $A\in\cL_Y\cap\cA(\Si)$ and find $\tA\in\cA_{\rm flat}(Y)$ such that 
$A=\tA|_\Si$.
Let $\a,\b\in\rT_A\cL_Y$, then by the characterization of $\cL_Y$ in (i) we 
find $\x,\z\in W^{1,p}(\Si,\cg)$ and paths 
$\tA^\a,\tA^\b :[-1,1] \to \cA_{\rm flat}(Y)$ such that 
$\tA^\a(0)=\tA^\b(0)=\tA$ and
$$
\a \;=\; 
\rd_A\x + \tfrac\rd{\ds}\bigr|_{\scriptscriptstyle s=0} \tA^\a (s) |_\Si ,
\qquad
\b \;=\; 
\rd_A\z + \tfrac\rd{\ds}\bigr|_{\scriptscriptstyle s=0} \tA^\b (s) |_\Si .
$$
Now firstly Stokes' theorem on $\Si$ with $\pd\S=\emptyset$ proves 
$$
\o(\rd_A\x\,,\,\rd_A\z)
\;=\; \lim_{\n\to\infty} \int_\Si \la \rd_A\x^\n \wedge \rd_A\z^\n \ra \\
\;=\; \lim_{\n\to\infty} \int_\Si \rd \la \x^\n \wedge \rd_A\z^\n \ra 
\;=\; 0.
$$
Here we have used smooth $W^{1,p}$-approximations $\x^\n$ and $\z^\n$ of $\x$ 
and $\z$ respectively.

Similarly, one obtains $\o(\rd_A\x,\frac\rd{\ds}\tA^\b|_\Si)=0$ and
$\o(\frac\rd{\ds}\tA^\a|_\Si,\rd_A\z)=0$ since
$\rd_A\bigl(\frac\rd{\ds}\tA^\a|_\Si\bigr)=\frac\rd{\ds}F_{\tA^\a}\bigr|_\Si=0$.
Finally, Stokes' theorem with $\pd Y=\Si$ yields due to $F_{\tA^\a(s)}=0$ for all $s$
\begin{align*}
\o(\a\,,\,\b) 
&\;=\; \int_\Si \la \tfrac\rd{\ds}\tA^\a|_\Si 
        \wedge \tfrac\rd{\ds}\tA^\b|_\Si \ra 
\;=\; \int_Y \rd \la \tfrac\rd{\ds}\tA^\a \wedge \tfrac\rd{\ds}\tA^\b \ra \\
&\;=\; \int_Y \la \tfrac\rd{\ds} F_{\tA^\a} \wedge \tfrac\rd{\ds} \tA^\b \ra 
      -\int_Y \la \tfrac\rd{\ds} \tA^\a \wedge \tfrac\rd{\ds} F_{\tA^\b} \ra 
\;=\,0 .
\end{align*}
This proves that $\o|_{\rT_A\cL_Y}=0$ and recalling (\ref{om,L2}) one moreover sees 
that $\rT_A\cL_Y$ and $*\rT_A\cL_Y$ are $L^2$-orthogonal. In fact, we even have
the topological decomposition 
$L^p(\Si,\rT^*\Si\otimes\cg) = \rT_A\cL_Y \oplus *\rT_A\cL_Y$,
and this proves the Lagrangian property of $\cL_Y$.
To see that this direct sum indeed exhausts the whole space consider the
Hodge type decomposition as in the proof of lemma~\ref{Lagrangian lemma},
$$
L^p(\Si,\rT^*\Si\otimes\cg) 
= \rd_A W^{1,p}_z(\Si,\cg) \oplus *\rd_A W^{1,p}_z(\Si,\cg) \oplus \tilde h_A .
$$
Here we have ${\rm dim}\,\tilde h_A = 2g \cdot {\rm dim}\,\rG$, and
we have already seen that $\cL_Y$ is a $\cG^{1,p}_z(\Si)$-bundle over the
$(g\cdot{\rm dim}\,\rG)$-dimensional manifold ${\rm Hom}(\pi_1(Y),\rG)$. So
$\rd_A W^{1,p}_z(\Si,\cg)\subset\rT_A\cL_Y$ is the tangent space to the fibre 
through $A$, and then for dimensional reasons $\rT_A\cL_Y \oplus *\rT_A\cL_Y$
also exhausts $\tilde h_A$ and thus all of $L^p(\Si,\rT^*\Si\otimes\cg)$.
\QED

 \bibliographystyle{alpha}

\begin{thebibliography}{14}
 \addcontentsline{toc}{section}{\numberline{}Bibliography}



 \bibitem[Ad]{Adams} R.A.Adams,
 {\it Sobolev Spaces}, Academic Press, 1978.

%
%

 \bibitem[At]{A1} M.F.~Atiyah,
 New invariants of three and four dimensional manifolds,
 {\it Proc. Symp. Pure Math.} {\bf 48} (1988).

 \bibitem[AB]{AB} M.F.Atiyah, R.Bott,
 The Yang Mills equations over Riemann surfaces,
 {\it Phil. Trans. R. Soc. Lond. A} {\bf 308} (1982), 523--615.

  \bibitem[B]{Broecker} T.Br\"ocker, T.tom~Dieck,
 {\it Representations of Compact Lie Groups}, Springer, 1985.

%
%
%
%
%
%
%

 \bibitem[DK]{DK} S.K.Donaldson, P.B.Kronheimer,
 {\it The Geometry of Four-Manifolds}, 
 Oxford Science Publications, 
 1990.

 \bibitem[DS1]{DS1} S.Dostoglou, D.A.Salamon,
 Instanton Homology and Symplectic Fixed Points,
 {Symplectic Geometry, edited by D. Salamon, Proceedings of a Conference}, 
 LMS Lecture Notes Series 192, Cambridge University Press, 1993, 57-94.

 \bibitem[DS2]{DS} S.Dostoglou, D.A.Salamon,
 Self-dual instantons and holomorphic curves,
 {\it Annals of Mathematics} {\bf 139} (1994), 581--640.

%
%
%

 \bibitem[F]{Fu} K.Fukaya,
 Floer homology for 3-manifolds with boundary I, Preprint 1997.
 \verb+http://www.kusm.kyoto-u.ac.jp/~fukaya/fukaya.html+

%
%
%
%
%
%
%
%
%
%

 \bibitem[L]{Lang} S.Lang,
 {\it Analysis II},
 Addison-Wesley, 1969.
 
 
 \bibitem[MS]{MS} D.McDuff, D.Salamon,
 {\it J-holomorphic Curves and Quantum Cohomology},
 new edition, in preparation.

 \bibitem[S]{Sa1} D.A.Salamon,
 Lagrangian intersections, $3$-manifolds with boundary, and the 
 Atiyah--Floer conjecture, 
 {\it Proceedings of the ICM, Z\"urich, 1994}, 
 Birkh\"auser, Basel, 1995, Vol.~1, 526--536.

%
%
%
%
%
%
%

 \bibitem[W1]{W} K.Wehrheim,
 {\it Uhlenbeck Compactness}, to appear in EMS Lectures in Mathematics. 

 \bibitem[W2]{W ell} K.Wehrheim,
 {\it Anti-self-dual instantons with Lagrangian boundary conditions I: Elliptic theory},
 preprint. 
 
%
%
 
 \end{thebibliography}

\end{document}